\documentclass[12pt]{amsart}
\usepackage{pdfsync}
\usepackage{graphicx}
\usepackage{amssymb,amsmath,latexsym,amsfonts,amsthm}
\usepackage{epstopdf}
\usepackage[mathscr]{eucal}
\DeclareMathOperator{\de}{d}

\DeclareMathOperator{\dist}{dist}

\begin{document}

\title{Oscillation of Urysohn type  spaces }
\author{ N. W. Sauer}

\date{January 04 2010 }
\maketitle

\newcommand{\snl} {\\ \smallskip}
\newcommand{\mnl}{\\ \medskip} 
\newcommand{\Bnl} {\\ Bigskip}
\newcommand{\edge}{\makebox[22pt]{$\circ\mspace{-6 mu}-\mspace{-6 mu}\circ$}}
\newcommand{\nedge}{\makebox[22pt]{$\circ\mspace{-6 mu}\quad\: \mspace{-6 mu}\circ$}}
\newcommand{\UR}{\boldsymbol{U}_{\hskip -3pt {\mathcal{R}}}}
\newcommand{\Ur}[1]{\boldsymbol{U}_{\hskip -3pt {#1}}}

\newcommand{\restrict}[2]{#1\mspace{-2mu}\mathbin{\upharpoonright}\mspace{-1mu} #2}
\newcommand{\Kat}{Kat\u{e}tov }
\newcommand{\Fra}{Fra\"{\i}ss\'e}
\newcommand{\str}[1]{\stackrel{#1}{\sim}}
\newcommand{\spe}{\mathrm{spec}}

\newcommand{\upl}[1]{\sideset{ }{^{c}}{\operatorname{\mathit{{#1}}}}}
\newcommand{\lpl}[1]{\sideset{ }{_{\negmedspace c}}{\operatorname{\mathit{{#1}}}}}
\newcommand{\ipl}[1]{\sideset{ }{_{\negmedspace c}^{c}}{\operatorname{\mathit{{#1}}}}}

\newcommand{\last}[1]{\mathfrak{#1}\!\rceil}

\newcommand{\bde}{\bar{\de}}
\newcommand{\dw}[1]{\ulcorner\negthickspace#1\negthickspace\urcorner}

\newtheorem{thm}{Theorem}[section]
\newtheorem*{thm*}{Theorem}
\newtheorem{lem}{Lemma}[section]
\newtheorem{coroll}{Corollary}[section]
\newtheorem{ass}{Assumption}[section]
\newtheorem{defin}{Definition}[section]
\newtheorem{example}{Example}[section] 
\newtheorem{fact}{Fact}[section]

\newtheorem{prop}{Proposition}[section]
\newtheorem{obs}{Observation}[section]
\newtheorem{cor}{Corollary}[section]
\newtheorem{sublem}{Sublemma}[section]
\newtheorem{claim}{Claim}
\newtheorem{question}{Question}[section]

\newtheorem{comment}{Comment}[section]
\newtheorem{problem}{Problem}[section]
\newtheorem{remark}{Remark}[section]

\newcommand\con{\char'136{}}

\begin{abstract}

A metric space $\mathrm{M}=(M;\de)$ is  {\em homogeneous} if for every isometry $\alpha$ of a finite subspace of\/ $\mathrm{M}$ to a subspace of\/ $\mathrm{M}$ there exists an isometry of $\mathrm{M}$ onto $\mathrm{M}$ extending $\alpha$. The metric space $\mathrm{M}$ is {\em universal} if it isometrically embeds every finite metric space $\mathrm{F}$ with $\dist(\mathrm{F})\subseteq \dist(\mathrm{M})$. ($\dist(\mathrm{M})$ being the set of distances between points of $\mathrm{M}$.)

A metric space $\mathrm{M}$ is {\em oscillation stable} if for every  $\epsilon>0$ and every  uniformly continuous and bounded function $f: M\to \Re$ there exists an isometric copy $\mathrm{M}^\ast=(M^\ast; \de)$ of\/ $\mathrm{M}$ in $\mathrm{M}$ for which: 
\[
\sup\{|f(x)-f(y)| \mid x,y\in M^\ast\}<\epsilon.
\]

\begin{thm*}
Every bounded, uncountable,  separable, complete, homogeneous, universal   metric space $\mathrm{M}=(M;\de)$   is oscillation stable.  (Theorem \ref{thm:finabstr}.)
\end{thm*}

\end{abstract}

\section{basic notions and facts}

Find in this section some basic definitions and facts and in the next, introductory  section, a general background discussion about the topics of this article.

\begin{defin}\label{defin:homuniv}
A metric space $\mathrm{M}=(M;\de)$ is  {\em homogeneous} if for every isometry $\alpha$ of a finite subspace of\/ $\mathrm{M}$ to a subspace of $\mathrm{M}$ there exists an isometry of $\mathrm{M}$ onto $\mathrm{M}$ extending $\alpha$. A metric space $\mathrm{M}$ is {\em universal} if it embeds every finite metric space $\mathrm{F}$ with $\dist(\mathrm{F})\subseteq \dist(\mathrm{M})$.
\end{defin}

\begin{defin}\label{defin:Urysohn}
A metric space $\mathrm{M}$ is an {\em Urysohn metric space}  if it is separable, complete, homogeneous and universal. 
\end{defin}

It follows easily from the general \Fra\  theory and is stated explicitly as Corollary 2.2 of \cite{SaCMJ} that 

\begin{fact}\label{fact:first}
An Urysohn metric space $\boldsymbol{U}$ embeds every separable metric space $\mathrm{N}$ with $\dist(\mathrm{N})\subseteq \dist(\boldsymbol{U})$. 
\end{fact}
See \cite{Ng} for further background information on homogeneous metric spaces and their partition theory. Two Urysohn metric spaces with the same set of distances are isometric. (Theorem \ref{thm:characterization}). An Urysohn metric space with $\mathcal{R}$ as set of distances will be denoted by $\UR$.

\begin{defin}\label{defin:osc}
A metric space $\mathrm{M}=(M;\de)$ is {\em oscillation stable} if for every  $\epsilon>0$ and every  uniformly continuous and bounded function $f: M\to \Re$ there exists a copy $\mathrm{M}^\ast=(M^\ast; \de)$ of\/ $\mathrm{M}$ in $\mathrm{M}$ for which:
\[
\sup\{|f(x)-f(y)| \mid x,y\in M^\ast\}<\epsilon.
\] 
\end{defin}

\begin{defin}\label{defin:appr}
Let $\mathrm{M}=(M;\de)$ be a metric space. Then for $\epsilon>0$ and $A\subseteq M$ let:
\begin{align}\label{align:epsnbhood}
\Bigl(A\Bigr)_\epsilon=\{x\in M\mid \exists y\in A\, \, (\de(x,y)<\epsilon)\}.
\end{align}
The metric space $\mathrm{M}$ is {\em approximately indivisible} if for every $\epsilon>0$ and $n\in \omega$ and function $\gamma: M\to n$ there exist $i\in n$ and an isometric copy $\mathrm{M}^\ast=(M^\ast;\de)$ of\/ $\mathrm{M}$ in $\mathrm{M}$ with
\[
M^\ast\subseteq \Bigl(\gamma^{-1}(i)\Bigr)_\epsilon.
\]
The metric space $\mathrm{M}$ is {\em  indivisible} if for every  $n\in \omega$ and function $\gamma: M\to n$ there exist $i\in n$ and an isometric copy $\mathrm{M}^\ast=(M^\ast;\de)$ of\/ $\mathrm{M}$ in $\mathrm{M}$ with
\[
M^\ast\subseteq \gamma^{-1}(i).
\]
\end{defin}

If $\mathcal{R}$ is finite then we obtain from Theorem 9.1 of \cite{Safinite}:
\begin{thm}\label{thm:finitedist}
Every Urysohn metric space $\UR$ for which the set of distances $\mathcal{R}$ is finite,  is indivisible. 
\end{thm}

\begin{thm*}{(Theorem \ref{thm:appequosc})}
A metric space $\mathrm{M}=(M;\de)$ is {\em oscillation stable} if and only if it is approximately indivisible. 
\end{thm*}
Hence we can use the terms oscillation stable and approximately indivisible interchangeably  to denote the same phenomenon. When citing theorems we will use original terminology. 

\begin{defin}\label{defin:oplus}
Let $\mathcal{R}\subseteq \Re$, then:
\[
a\oplus b:=\sup\{x\in \mathcal{R} \mid x\leq a+b\}.
\] 
\end{defin}

\section{Introduction}

Sharpening and reformulating Dvoretzky's Theorem, see \cite{Dvoretzky1} and \cite{Dvoretzky2}, V. Milman, see \cite{Mil}, proved, for  $\mathbb{S}^k$ the unit sphere of the $k+1$-dimensional Euclidean space and for $\mathbb{S}^\infty$, the unit sphere of the Hilbert space $\ell_2$:

\begin{thm*}[Milman]
\label{thm:Milman'} Let $n\in \omega$ and  $\gamma:\mathbb{S} ^{\infty}\to n$. Then for every
$\varepsilon > 0$ and every  $k \in \omega$, there is $i\in n$ and an isometric copy
$(\mathbb{S}^k)^\ast$ of $\mathbb{S}^k$ in $\mathbb{S} ^{\infty}$ such that
$(\mathbb{S}^k)^\ast \subseteq \Bigl(\gamma^{-1}(i)\Bigr)_\epsilon$.
\end{thm*}

This result then led naturally to the {\em distortion problem}, asking if the Hilbert sphere $\mathbb{S}^\infty$ is approximately indivisible, that is oscillation stable. Which was settled by:

\begin{thm*}[Odell-Schlumprecht \cite{OS}]

\label{thm:Odell-Schlumprecht}

There is a number $n\in \omega$ and a function $\gamma: \mathbb{S}^\infty\to n$   and  $\epsilon > 0$ such for every $i\in n$ the set $\Bigl(\gamma^{-1}(i)\Bigr)_\epsilon$ does not contain a copy of $\mathbb{S} ^{\infty}$. (That is $\mathbb{S}^\infty$ has {\em distortion}.)
\end{thm*}

The notion of oscillation stable topological groups was first introduced and studied by Kechris, Pestov and Todorcevic in \cite{KPT}, see also Pestov \cite{Pe1}. For  homogeneous metric spaces, their groups of isometries are oscillation stable if the metric spaces are oscillation stable in the sense of Definition  \ref{defin:osc}. In the special case of metric spaces it is more convenient to use the notion of oscillation stable as given in  Definition~\ref{defin:osc}. It is shown in \cite{KPT}, see also \cite{Pe1}, that homogeneous metric spaces are oscillation stable if and only if they are approximately indivisible. It follows from Theorem \ref{thm:appequosc} that this   equivalence actually holds for all metric spaces. 

The Hilbert Sphere $\mathbb{S}^\infty$ and the Urysohn Sphere $\Ur{[0,1]}$ share many interesting topological properties. For example:  The group $\mathrm{Iso}(\mathbb{S}^\infty)$ is extremely amenable, as shown by Gromov and Milman \cite{GM}. The group $\mathrm{Iso}(\Ur{[0,1]}$ is extremely amenable, as shown by Pestov \cite{Pe0}. Subsequently to the Odell-Schlumprecht result, it was therefore natural to ask whether the Urysohn sphere $\Ur{[0,1]}$ has distortion as well. The first major step in resolving this question is the main result achieved by  Lopez-Abad and  Nguyen Van Th\'e in \cite{LANVT}: The Urysohn sphere is oscillation stable if all Urysohn metric spaces $\Ur{n}$ for $n=\{0,1,\dots,n-1\}$ are indivisible. This then was proven in \cite{LNSA} by Nguyen Van Th\'e  and Sauer, yielding the result that the Urysohn sphere $\Ur{[0,1]}$ is  oscillation stable. 

The main result in this paper, namely that bounded, uncountable Urysohn metric spaces are oscillation stable, required  a sequence of results similar to the  Abad and  Nguyen Van Th\'e and Sauer sequence of results. 

The characterization of the distance sets of Urysohn metric spaces was achieved in \cite{SaCMJ}, see Section \ref{sect:4-values} Theorem \ref{thm:characterization} of this article. In particular the set of distances $\mathcal{R}$ of an Urysohn metric space is closed and is closed under the operation $\oplus$ which is then associative on $\mathcal{R}$. 

Let $\mathcal{R}$ be a closed set of numbers for which $\oplus$ is associative. As a first step we need to approximate $\mathcal{R}$ with a finite subset for which $\oplus$ is associative. This is accomplished in Section \ref{sect:finappr} Thorem \ref{thm:fext}. It follows from Theorem 9.1 of \cite{Safinite} that  Urysohn metric spaces with a finite set of distances are indivisible, extending the result that the Urysohn metric spaces $\Ur{n}$ are indivisible.   

An essential part of the proof in the Lopez-Abad--Nguyen paper is a clever tree construction. For this construction it was necessary to calculate certain distances.  This was possible because the finite approximating sets of numbers for the interval $[0,1]$ are sets of the form $\{\frac{m}{n}\mid 0\leq m\leq n\}$, and hence explicitly available for the calculation. In the general case dealt with in this article, the finite approximating sets are not explicitly available and hence it is impossible to calculate the required distances.  It follows that the existence of such distances has to be proven. This necessitated the study of metric $\mathcal{R}$-graphs in Section \ref{sect:metricgra}. 

Finally we have to address the question: What are the  closed subsets of the reals having 0 as a limit and for which $\oplus$ is associative? Some initial experiments seemed to indicate that such sets are just finite unions of closed intervals. If this where  the case, then  the oscillation theorem for Urysohn metric spaces could simply be derived from the  Abad--Nguyen Van Th\'e--Sauer result.  But, as will be shown in Section~\ref{sect:Cantor-4val}, there are quite intricate examples,  of such sets. (Generalized Cantor type sets for example.) A complete characterization of such sets seems to be beyond our present abilities.

\section{The 4-values condition}\label{sect:4-valuescond}\label{sect:4-values}

Let $\mathcal{R}\subseteq \Re_{\geq 0}$ be a closed subset of the non negative reals.  Then, according to Definition \ref{defin:oplus}, the operation $\oplus$ or   $\oplus_\mathcal{R}$ if the distinction is needed is given by: 
\[
a\oplus b:=\sup\{x\in \mathcal{R} \mid x\leq a+b\}
\] 
which, because $\mathcal{R}$ is closed as a subset of $\Re$, is a binary operation on $\mathcal{R}$.    A subset $\{a,b,c\}\subseteq \mathcal{R}$ is {\em metric} if $a\leq b+c$ and $b\leq a+c$ and $c\leq a+b$. A triple $(a,b,c)$ of numbers in $\mathcal{R}$ is {\em metric} if its set of entries $\{a,b,c\}$ is metric.   Note the following immediate consequence:

\begin{obs}\label{obs:opl}
$\max\{b,c\}\leq b\oplus  c\in \mathcal{R}$ and the triple $(b\oplus c, b,c)$ is metric, for all $b,c\in \mathcal{R}$. Also: 

If  $(a,b,c)$ is  a metric triple of numbers in $\mathcal{R}$, then $b\oplus c\geq a$ and   hence  $b\oplus c$ is the largest number  $x\in\mathcal{R}$ for which the triple $(x,b,c)$ is metric. The operation $\oplus$ is commutative and if $a\geq b$ then $a\oplus c\geq b\oplus c$. 
\end{obs}

\begin{defin}\label{defin:4value}
The set of quadruples with entries in $\mathcal{R}$ and with $\max\{b,c,d\}\leq a\leq b+c+d$ will be denoted by $\mathscr{Q}(\mathcal{R})$. The arrow $x\leadsto (a,b,c,d)$ means that the triples $(a,b,x)$ and $(c,d,x)$ are metric. The set $\mathcal{R}$ satisfies the {\em 4-values condition} if:
\begin{align*}
&\text{For all quadruples $(a,b,c,d)\in \mathscr{Q}(\mathcal{R})$ and all $x\in \boldsymbol{R}$}\\
&x\leadsto (a,b,c,d) \text{ \  implies that there exists $y\in \mathcal{R}$ with $y\leadsto (a,d,c,b)$}.
\end{align*}
\end{defin} 
The importance of the 4-values condition when investigating  homogeneous  metric spaces was discovered in \cite{DLPS} and then used in \cite{Ng} and \cite{SaCMJ}. The connection between the operation $\oplus$ being associative and the 4-values condition stated in the next  Theorem \ref{lthm:associative}, together with the result on $\mathcal{R}$-graphs stated in Theorem \ref{thm:graphext}, allowes different and sometimes simpler arguments for the results in \cite{SaCMJ}.

\begin{thm}\label{lthm:associative}
A closed subset $\mathcal{R}\subseteq \Re_{\geq 0}$  of the non negative reals satisfies  the the 4-values condition  if and only if the operation $\oplus$ on $\mathcal{R}$ is associative.   
\end{thm}
\begin{proof}
Let $\mathcal{R}$ satisfy the 4-values condition  and  $\{b,c,d\}\subseteq \mathcal{R}$.   We will prove that $d\oplus(c\oplus d)=(d\oplus c)\oplus b$. 

It follows from Observation \ref{obs:opl} that $d\oplus (c\oplus  b)=(b\oplus c)\oplus d$, that the triples $(d\oplus c, d,c)$ and  $\big((d\oplus c)\oplus b, d\oplus c, b\big)$ are metric and that $d+c+b\geq (d\oplus c)\oplus b\geq \max\{d,c,b\}$. Hence $\big((d\oplus c)\oplus b, b,c,d\big)\in \mathscr{Q}(\mathcal{R})$ and $d\oplus c\leadsto \big((d\oplus c)\oplus b, b,c,d\big)$. 

Because $\mathcal{R}$ satisfies the 4-values condition there exists $y\in \mathcal{R}$ with  $y\leadsto \big((d\oplus c)\oplus b, d, c,b\big)$.  Then $y\leq b\oplus c$ because the triple $(y,b,c)$ is metric and $y\oplus d\geq (d\oplus c)\oplus b$ because the triple $\big((d\oplus c)\oplus b, y,d\big)$ is metric. Hence:
\[
d\oplus (c\oplus  b)=(b\oplus c)\oplus d\geq y\oplus d\geq (d\oplus c)\oplus b,   
\]
and therefore:
\[
(d\oplus c)\oplus b=b\oplus(c\oplus d)\geq (b\oplus c)\oplus d=d\oplus (c\oplus b) \geq (d\oplus c)\oplus b,
\]
implying $d\oplus(c\oplus d)=(d\oplus c)\oplus b$.

Let $\oplus$ be associative and $(a,b,c,d)\in \mathscr{Q}(\mathcal{R})$ with $x\in \mathcal{R}$ and $x\leadsto (a,b,c,d)$. Then $d\oplus c\geq x$ and $x\oplus b\geq a$ and 
\[
(b\oplus c)\oplus d=d\oplus(c\oplus b)=(d\oplus c)\oplus b\geq x\oplus b\geq a.
\]
It follows that $b\oplus c\leadsto (a,d,c,b)$. 
\end{proof}

Theorem \ref{lthm:associative} together with 
\begin{thm*}( \cite{SaCMJ} Theorem 2.2.)\\
Any two Urysohn metric spaces $\mathrm{M}$ and $\mathrm{N}$  with $\dist(\mathrm{M}) = \dist(\mathrm{N})$ are isometric.
\end{thm*}
and
\begin{thm*}{( \cite{SaCMJ}Theorem 4.4.)}\\
Let $0\in \mathcal{R}\subseteq \Re_{\geq 0}$ with 0 as a limit. Then there exists an Urysohn metric space $\UR$ if and only if $\mathcal{R}$ is a closed subset of $\Re_{\geq 0}$ which satisfies the 4-values condition.

Let $0\in \mathcal{R}\subseteq \Re_{\geq 0}$ which does not have 0 as a limit. Then there exists an Urysohn metric space $\UR$  if and only if $\mathcal{R}$ is a countable subset of $\Re$ which satisfies the 4-values condition.
\end{thm*}
imply:
\begin{thm}\label{thm:characterization}
Any two Urysohn metric spaces $\mathrm{M}$ and $\mathrm{N}$ with $\dist(\mathrm{M})=\dist(\mathrm{N})$ are isometric.

Let $0\in \mathcal{R}\subseteq \Re_{\geq 0}$ with $0$  a limit of $\mathcal{R}$.  Then there exists an Urysohn metric space $\UR$ with $\dist(\UR)=\mathcal{R}$ if and only if $\mathcal{R}$ is   a closed subset of\/ $\Re_{\geq 0}$ with $\oplus_\mathcal{R}$ associative. 

Let $0\in \mathcal{R}\subseteq \Re_{\geq 0}$ which does not have 0 as a limit. If $\mathcal{R}$ is closed then there exists an Urysohn metric space $\UR$  if and only if $\mathcal{R}$ is a countable subset of $\Re$ with $\oplus_\mathcal{R}$ associative. 
\end{thm}
An Urysohn metric space with $\mathcal{R}$ as set of distances will be denoted by $\UR$.

\section{Ordered Urysohn metric spaces}\label{sect:orderedsp}

The purpose of this section is to establish Lemma \ref{lem:ordemb}, which will be needed later on. For readers familiar with the general \Fra\  theory the Lemma is probably quite obvious, but unfortunately it does not seem to be stated explicitly in the literature. Using  the tools of \Fra \  theory one would use the fact that  the class of finite induced substructures of a homogeneous structure has amalgamation and then conclude from Proposition~1.3 of \cite{DLPS} or Theorem 5 of \cite{Ng} that  in the case of metric spaces the amalgamation property implies the disjoint amalgamation property. It is then easily seen that the disjoint amalgamation property implies that every \Kat functions has infinitely many realizations. This then in turn implies Lemma \ref{lem:ordemb} via a standard construction  as for the special case proven in \cite{LANVT}.

Avoiding \Fra \ theory we use Theorem 2 of \cite{ELSaRamsey} which says that for every homogeneous structure for which the class of finite induced substructures has disjoint amalgamation, hence in particular for every Urysohn metric space $\boldsymbol{U}=(U;\de)$ and every finite subset $F\subseteq U$, the restriction $\restrict{\boldsymbol{U}}{(U\setminus F)}$ is an isometric copy of $\boldsymbol{U}$. A relational structure having this property is called {\em strongly inexhaustible}. Hence every Urysohn metric space is   strongly inexhaustible.

\begin{defin}\label{defin:omorder}
Let $\boldsymbol{U}=(U;\de)$ be a countable Urysohn metric space with $U=\{u_i\mid i\in \omega\}$ an enumeration of $U$. The {\em enumeration-order} on $U$ is the linear order given by: $u_i\leq u_j$ if and only if $i\leq j$. 
\end{defin}

\begin{lem}\label{lem:ordemb}
Let $\boldsymbol{U}=(U;\de)$ be a countable Urysohn metric space with $U=\{u_i\mid i\in \omega\}$ an enumeration of $U$ and let $\boldsymbol{U}^\ast=(U^\ast;\de)$ be an isometric copy of $\boldsymbol{U}$ in $\boldsymbol{U}$. Then there exists an enumeration-order preserving isometry of $\boldsymbol{U}$ into $\boldsymbol{U}^\ast$. 
\end{lem}
\begin{proof}
For $\alpha_0$ the empty isometry we will construct an $\omega$-sequence of isometries $\alpha_0\subseteq \alpha_1\subseteq \alpha_2\subseteq \alpha_3\subseteq$ with $\alpha_n: \{u_i\mid i\in n\}\to U^\ast$  which is enumeration-order  preserving for each $n\in \omega$. Then $\alpha:=\bigcup_{n\in \omega}\alpha_n$ will be an order preserving isometry of $\boldsymbol{U}$ into $\boldsymbol{U}^\ast$.

For $\alpha_n$ constructed let $l\in \omega$ such that $\alpha(n-1)=u_l$ and let 
\[
F=\big\{u_j\in U^\ast\mid j<l \text{  and $u_j\not\in \{\alpha(u_i)\mid i\in n\}$}\big\}
\]
and let $\bar{\boldsymbol{U}}=\restrict{\boldsymbol{U}^\ast}{(U^\ast\setminus F)}$. Let $\mathrm{N}$ be the subspace of $\boldsymbol{U}$ induced by the set $\{u_i\mid 0\leq i\leq n\}$. Because $\boldsymbol{U}$ is strongly inexhaustible, the space $\bar{\boldsymbol{U}}$ is isometric to $\boldsymbol{U}$ and hence universal and therefore there exists an isometry $\beta$ of the space $\mathrm{N}$ into $\bar{\boldsymbol{U}}$. Let $\gamma$ be the restriction of $\beta$ to the set $\{u_i\mid i\in n\}$. Then $\gamma^{-1}\circ \alpha_n$ is an isometry and because $\bar{\boldsymbol{U}}$ is homogeneous it has an extension $\delta$ to an isometry of $\bar{\boldsymbol{U}}$ onto $\bar{\boldsymbol{U}}$. Let  $\alpha_{n+1}=\delta\circ\beta$. 
\end{proof}

\section{Metric $\mathcal{R}$-graphs}\label{sect:graphs}\label{sect:metricgra}

Let $\mathcal{R}\subseteq \Re_{\geq 0}$ be a closed subset of the non negative reals satisfying  the 4-values condition. Then the $\oplus $ operation on $\mathcal{R}$ is commutative and  associative according to Theorem \ref{lthm:associative}. For $(a_i\in \mathcal{R}; i\in n\in \omega)$ let $\bigoplus_{i\in n}a_i=a_0\oplus a_1\oplus a_2\oplus \dots \oplus a_{n-1}$. 

\begin{lem}\label{lem:infoplus}
Let $a\in \mathcal{R}$ and $\mathcal{S}\subseteq \mathcal{R}$. Then $a\oplus\inf(\mathcal{S})=\inf\{a\oplus s\mid s\in \mathcal{S}\}$.
\end{lem}
\begin{proof}
If $a+\inf(\mathcal{S})=a\oplus\inf(\mathcal{S})$ then $a+s\geq a\oplus s\geq a\oplus\inf(\mathcal{S})=a+\inf(\mathcal{S})$ for all $s\in \mathcal{S}$ and hence $a+\inf(\mathcal{S})=\inf\{a+ s\mid s\in \mathcal{S}\}\geq \inf\{a\oplus s\mid s\in \mathcal{S}\}\geq a\oplus\inf(\mathcal{S})=a+\inf(\mathcal{S})$ implying  $a\oplus\inf(\mathcal{S})=\inf\{a\oplus s\mid s\in \mathcal{S}\}$.

Otherwise $a+\inf(\mathcal{S})>a\oplus\inf(\mathcal{S})$. Because $\mathcal{R}$ is closed there exists an $\epsilon>0$ with $[a\oplus \inf(\mathcal{S}), a+\inf(\mathcal{S})+\epsilon]\cap \mathcal{R}=\emptyset$. Then $a\oplus s=a\oplus\inf(\mathcal{S})$ for all $s\in \mathcal{S}$ with $s-\inf(\mathcal{S})<\epsilon$. Hence $\inf\{a\oplus s\mid s\in \mathcal{S}\}=a\oplus\inf(\mathcal{S})$.
\end{proof}

A sequence $\mathrm{V}=(v_0,v_1,v_2,\dots,v_{n-1},v_n)$ of vertices of a graph is a {\em walk from $v_0$ to $v_n$} if $\{x_i,x_{i+1}\}$ is an edge of the graph for all $i\in n$. The walk $P$ is a {\em trail} from $v_0$ to $v_n$ if $v_i=v_j$ only for $i=j$. A {\em path} $\{v_0, v_1, \dots, v_n\}$ of a graph is a trail in which two vertices $x_i$ and $x_j$  are adjacent if and only if $|i-j|=1$. Note that every walk from $a$ to $b$ contains a path from $a$ to $b$. The {\em length} of a path is the number of its edges.  A graph is {\em connected} if for all vertices $a$ and $b$ there is a path from $a$ to $b$.    A {\em cycle} is a connected graph in which every vertex is adjacent to exactly two points. (Hence, for our purposes,  a single edge is not a cycle of length two.)  It follows that the number of edges of a cycle is equal to the number of its vertices, the {\em length} of the cycle.  If $\{a,b\}$ is an edge of a cycle  of length $n$ then, after removing the edge $\{a,b\}$ from the cycle,  there is a path from $a$ to $b$  of length $n-1$.  An isomorphic embedding of a cycle into a graph is a {\em  cycle of the graph.} That is, a cycle of a graph  does not have any secants. Note that the vertices of a cycle of a graph can be adjacent to more than two points, just not in the cycle.

\begin{defin}\label{defin:Rgraph}
An {\em $\mathcal{R}$-graph} $\mathrm{G}=(G; E, \de)$ is a triple, for which:
\begin{enumerate}
\item $(G;E)$ is a simple graph with $G$ as set of vertices and $E$ as set of edges.
\item  $\de: E\cup \{(x,x)\mid x\in G\}\to \mathcal{R}$, the {\em distance function of\/ $\mathrm{G}$},   is a function with $\de(x,y)=0$ if and only if $x=y$.
\item $\de(x,y)=\de(y,x)$ for all $\{x,y\}\in E$. 
\end{enumerate}
The $\mathcal{R}$-graph   $\mathrm{G}=(G; E, \de)$ is {\em connected} if the graph $(G;E)$ is connected. 
\end{defin}
Note that if $\mathrm{M}=(M;\de)$ is a metric space with $\dist(\mathrm{M})\subseteq \mathcal{R}$ then $\mathrm{G}=(M;[M]^2,\de)$ is an $\mathcal{R}$-graph for which the underlying  graph $(M;[M]^2)$ is complete.

For $\mathrm{G}=(G; E, \de)$ an $\mathcal{R}$-graph and  $\mathrm{W}=(v_0,v_1,v_2,\dots,v_{n-1},v_n)$ a walk let $\de(W)=\bigoplus_{i\in n} \de(v_i,v_{i+1})$.  Note that for every walk $W$ from $v_0$ to $v_n$ there exists a trail $P$ from $v_0$ to $v_n$ with $\de(P)\leq \de(W)$.  
For $\mathrm{G}=(G;E, \de)$ and  $a,b\in G$  let: 
\begin{align*}
&\boldsymbol{d}(a,b):=\inf\{\de(\mathrm{W})\mid \text{ $\mathrm{W}$ is a walk from $a$ to $b$}\}.
\end{align*} 
Note that $\boldsymbol{d}(a,b)=\inf\{\de(\mathrm{P})\mid \text{ $\mathrm{P}$ is a trail from $a$ to $b$}\}$.

\begin{defin}\label{defin:metricdgr}
An $\mathcal{R}$-graph $\mathrm{G}=(G; E, \de)$ is  {\em regular} if $\boldsymbol{d}(a,b)>0$ for all $(a,b)\in G^2$ with $a\not=b$. (Note that   finite $\mathcal{R}$-graphs are regular.)

An $\mathcal{R}$-graph $\mathrm{G}=(G; E, \de)$ is  {\em metric} if $\de(v_0,v_{n})\leq\de(W)$ for every walk  $(v_0,v_1,\dots,v_{n-1},v_n)$  from $v_0$ to $v_n$ and  every edge $\{v_0,v_n\}\in E$.  That is if $\boldsymbol{d}(a,b)=\de(a,b)$  for every edge $\{a,b\}\in E$.
\end{defin}

It follows that a cycle $\mathrm{C}=(C;E_\mathrm{C}, \de)$  is metric if and only if  for every edge $e\in E_{\mathrm{C}}$: 
\[
\de(e)\leq \bigoplus_{h\in E_\mathrm{C}\setminus\{e\}} \de(h).
\] 

\begin{lem}\label{lem:cyclemetric}
 A regular $\mathcal{R}$-graph $\mathrm{G}=(G; E, \de)$ is metric if and only if all of its  cycles are metric. 
\end{lem}
\begin{proof}
If $\mathrm{G}$ contains a cycle induced by the set $\{v_i\mid i\in n\in \omega\}$ for which $\de(v_0,v_n)>\bigoplus_{i\in n}\de(v_i,v_{i+1})$ then $\mathrm{G}$ is not metric because $v_0,v_1,\dots, v_n$ is a walk from $v_0$ to $v_n$. Assume that all cycles of $\mathrm{R}$ are metric and for a contradiction that the $\mathcal{R}$-graph $\mathrm{G}$ is not metric. 

Let $n\in \omega$ be the smallest number for which there exists an edge $\{v_0,v_1\}\in E$ and a trail $\mathrm{T}=(v_0,v_1,\dots, v_n)$ with $\de(v_0,v_1)>\de(T)$. Then $\{v_i\mid i\in n\}$ does not induce a cycle and hence there exists a secant, that is  there are indices $i,j$ with $\{i,j\}\not=\{0,n\}$ and with  $i+1<j$ so that $\{v_i,v_j\}$ is an edge in $E$.  Then 
\[\de(v_i,v_j)\leq \de(v_i,v_{i+1})\oplus\de(v_{i+1},v_{i+2})\oplus\dots\oplus\de(v_{j-1},v_j)
\]
follows from  $j-i<n$ and the minimality condition on $n$. The trail $\mathrm{S}=(v_0,v_1,\dots, v_i,v_j,v_{j+1},\dots,v_n)$ is shorter than the trail $\mathrm{T}$ with $\de(S)\leq \de(T)<\de(v_0,v_1)$ in contradiction to the minimality of $n$. 
\end{proof}

\begin{cor}\label{cor;mspmetric}
Let $\mathrm{M}=(M;\de)$ be a metric space with $\dist(\mathrm{M})\subseteq \mathcal{R}$. Then the $\mathcal{R}$-graph $\mathrm{G}=(M; [M]^2,\de)$ is metric and regular.
\end{cor}
\begin{proof}
The only cycles of $\mathrm{G}$ are triangles, which are metric. Hence it follows from Lemma \ref{lem:cyclemetric} that the $\mathcal{R}$-graph $\mathrm{G}$ is metric and regular. 
\end{proof}

For two $\mathcal{R}$-graphs $\mathrm{G}=(G;E_\mathrm{G}, \de^\mathrm{G})$ and $\mathrm{H}=(G;E_\mathrm{H},\de^\mathrm{H})$ let $\mathrm{G}\curlyeqprec \mathrm{H}$ if $E_\mathrm{G}\subseteq E_\mathrm{H}$ and $\de^\mathrm{H}$ restricted to $E_\mathrm{G}$ is equal to $\de^\mathrm{G}$.    If $\mathrm{G}=(G; E,\de^\mathrm{G})$ is an $\mathrm{R}$-graph and $\mathrm{M}=(G;\de^\mathrm{M})$ a metric space with $\dist(\mathrm{M})\subseteq \mathcal{R}$ and $\mathrm{H}$ the $\mathcal{R}$-graph $(G;[G]^2, \de^\mathrm{M})$  then we write $\mathrm{G}\curlyeqprec \mathrm{M}$ if $\mathrm{G}\curlyeqprec \mathrm{H}$. Note that if $\mathrm{G}\curlyeqprec \mathrm{H}$  and $\mathrm{H}$ is metric and regular then $\mathrm{G}$ is metric and regular.

\begin{lem}\label{lem:compdist}
Let $\mathrm{G}=(G; E, \de^{\mathrm{G}})$ be a connected and metric regular $\mathcal{R}$-graph with $a,b\in G$ and $a\not= b$ and $\{a,b\}\not\in E$.

Then  $\mathrm{H}=(G; E\cup \{a,b\}, \de^{\mathrm{H}})$ with $\mathrm{G}\curlyeqprec \mathrm{H}$ and $\de^{\mathrm{H}}(a,b)=\boldsymbol{d}^{\mathrm{G}}(a,b)$ is a metric regular $\mathcal{R}$-graph with $\boldsymbol{d}^\mathrm{G}(x,y)=\boldsymbol{d}^\mathrm{H}(x,y)$  for all $\{x,y\}\in [G]^2$.
\end{lem}
\begin{proof}
Let $\{x,y\}\in [G]^2$. In order to determine  $\boldsymbol{d}^\mathrm{H}(x,y)$ let   $\mathcal{P}$ be the set of trails in $\mathrm{H}$  from $x$ to $y$.  Let $\mathcal{Q}$ be the set of trails from $x$ to $y$  in $\mathrm{G}$ and let $\mathcal{Z}$ be the set of trails in $\mathrm{G}$ from $a$ to $b$.  Let $\mathrm{P}=(x=v_0,v_1,v_2,\dots,v_{n-1},v_n=y)\in \mathcal{P}$. If $\{a,b\}$ is not an edge of $\mathrm{P}$ then $\de^{\mathrm{H}}(\mathrm{P})=\de^{\mathrm{G}}(\mathrm{P})$ and $\mathrm{P}\in \mathcal{Q}$.   Otherwise there is exactly one $i\in n$ so that $\{v_i,v_{i+1}\}=\{a,b\}$. Let $\mathrm{T}$ be the trail $v_0,v_1,\dots, v_i$ in $\mathrm{G}$ and $\mathrm{S}$ the trail $v_{i+1}, v_{i+2}, \dots, v_n$ in $\mathrm{G}$ and let $a=\de^\mathrm{G}(\mathrm{T})+\de^\mathrm{G}(\mathrm{S})=\de^\mathrm{H}(\mathrm{T})+\de^\mathrm{H}(\mathrm{S})$.

Then, using Lemma \ref{lem:infoplus} 
\[
\de^\mathrm{H}(\mathrm{P})=a\oplus \boldsymbol{d}^\mathrm{G}(a,b)=\inf\{a\oplus\de^{\mathrm{G}}(\mathrm{Z})\mid \mathrm{Z}\in \mathcal{Z}\}.
\] 
Hence there exists, for every $\epsilon>0$, a trail $\mathrm{Z}\in \mathcal{Z}$ so that $\de^\mathrm{G}(\mathrm{Q})\leq \de^\mathrm{H}(\mathrm{P})+\epsilon$ for $\mathrm{Q}$ the trail in $\mathrm{G}$ from $x=v_0$ to $y=v_n$ which goes from $v_0$ to $v_i$ along $\mathrm{T}$ and then to $v_{i+1}$ along $\mathrm{Z}$ and then to $v_n=y$ along $\mathrm{S}$. Hence, for every $\epsilon>0$ and for every trail $\mathrm{P}\in \mathcal{P}$ there exists a trail $\mathrm{Q}_\mathrm{P}\in \mathcal{Q}$  with $\de^{\mathrm{G}}(\mathrm{Q}_\mathrm{P})<\de^{\mathrm{H}}(\mathrm{P})+\epsilon$. This together with $\mathcal{Q}\subseteq \mathcal{P}$ implies that 
\[
\boldsymbol{d}_\mathrm{H}(x,y)=\inf\{\de^{\mathrm{H}}(\mathrm{P})\mid \mathrm{P}\in \mathcal{P}\}=\inf\{\de^{\mathrm{G}}(\mathrm{Q})\mid \mathrm{Q}\in \mathcal{Q}\}=\boldsymbol{d}_\mathrm{G}(x,y).
\]  
It follows that $\mathrm{H}$ is metric and regular.

\end{proof}

\begin{lem}\label{lem:connext}
Let  $\mathrm{G}=(G; E_\mathrm{G}, \de^\mathrm{G})$ be a metric regular $\mathcal{R}$-graph. There exists a connected,  metric regular $\mathcal{R}$-graph $\mathrm{H}=(G; E_\mathrm{H},\de^\mathrm{H})$ with $\mathrm{G}\curlyeqprec\mathrm{H}$.
\end{lem}
\begin{proof}
Let $\mathfrak{C}$ be the set of connected components of $\mathrm{G}$.   For each $\mathrm{C}\in \mathfrak{C}$ let $v_\mathrm{C}$ be a vertex in the connected component $\mathrm{C}$.  Let $0<r\in \mathcal{R}$ and $V=\{v_\mathrm{C}\mid \mathrm{C}\in \mathfrak{C}\}$. Let  $\mathrm{H}=(G; E_\mathrm{H}, \de^\mathrm{H})$ the $\mathcal{R}$-graph with $E_\mathrm{H}=E_\mathrm{G}\cup [V]^2$ and $\de(v_\mathrm{C},v_\mathrm{D})=r$ for all edges $\{v_\mathrm{C},v_\mathrm{D}\}\in [V]^2$ and with $\mathrm{G}\curlyeqprec\mathrm{H}$.

A trail which contains points in different connected components contains an edge of distance $r>0$. Hence $\mathrm{H}$ is regular because each of the connected components is reguar. If a cycle of $\mathrm{H}$ contains two vertices in $V$ then it is a subset of $V$ and hence a metric triangle. Otherwise the cycle is a subset of one of the connected components and by assumption metric. 

\end{proof}

\begin{lem}\label{lem:graphext}
Let $\mathrm{G}=(G; E, \de)$ be a connected, regular,  metric  $\mathcal{R}$-graph. Then $\mathrm{M}=(G,\boldsymbol{d}^\mathrm{G})$ is a metric space on $G$ with $\de^\mathrm{G}(a,b)=\boldsymbol{d}^\mathrm{G}(a,b)$ for all edges $\{a,b\}\in E(\mathrm{G})$.
\end{lem}
\begin{proof}
Let $\mathcal{P}$ be the partial $\curlyeqprec$-order of metric and regular $\curlyeqprec$-extensions $\mathrm{H}$  of the $\mathcal{R}$-graph $\mathrm{G}=(G; E, \de)$ with $\boldsymbol{d}^\mathrm{H}(a,b)=\boldsymbol{d}^\mathrm{G}(a,b)$ for all points $a,b\in G$. Let $\mathcal{C}=\big(\mathrm{H}_i = (G; E_i, \de^i); i\in I\big)$ be a chain of $\mathcal{P}$ and $\mathrm{K}=(G; E_\mathrm{K},\de^\mathrm{K})$ the $\mathcal{R}$-graph with $E_\mathrm{K}=\bigcup_{i\in I}E_i$ and $\de^\mathrm{K}(x,y)=\de^i(x,y)$ for some $i\in I$ with $\{x,y\}\in E_i$.  Let $\mathrm{V}$ be a cycle of $\mathrm{K}$. Because $\mathrm{V}$ is finite there is an $i\in I$ so that $\mathrm{V}$ is a cycle of $\mathrm{H}_i$ and hence is metric. It follows now from Lemma \ref{lem:cyclemetric} that  $\mathrm{K}$ is metric.

We have to check that $\mathrm{K}$ is regular. Assume for a contradiction that there are two vertices $a,b\in G$ with $a\not=b$ and a sequence $(T_n)$ of trails in $\mathrm{K}$  from $a$ to $b$ so that the sequence $\de^\mathrm{K}(T_n)$ tends to 0. For each of the trails $T_n$ exists an $i\in I$ so that $T_n$ is a trail in $\mathrm{H}_i$ and hence $\de^\mathrm{K}(T_n)=\de^{\mathrm{H}_i}(T_n)\geq \boldsymbol{d}^{\mathrm{H}_i}(a,b)=\boldsymbol{d}^\mathrm{G}(a,b)$, a contradiction. 

A simple application of Zorn's Lemma to the partial order $\mathcal{P}$  together with Lemma~\ref{lem:compdist} shows that the metric $\mathcal{R}$-graph $\mathrm{G}$ has a $\curlyeqprec$-extension to a metric space  $\mathrm{M}$ on $G$. 
\end{proof}

\begin{thm}\label{thm:graphext}
Let $\mathcal{R}\subseteq \Re_{\geq 0}$ be a closed subset of the non negative reals satisfying  the 4-values condition and let $\mathrm{G}=(G; E_\mathrm{G}, \de^\mathrm{G})$ be an  $\mathcal{R}$-graph. 

If there exists a metric space $\mathrm{M}=(G;\de^\mathrm{M})$ with $\de^\mathrm{G}(a,b)= \de^\mathrm{M}(a,b)$ for all edges $\{a,b\}\in E_\mathrm{G}$, then the $\mathcal{R}$-graph $\mathrm{G}$ is  regular and  metric.

If $\mathrm{G}$ is metric and regular, then there exists a connected, metric and regular $\mathcal{R}$-graph $\mathrm{H}=(G;E_\mathrm{H}, \de^\mathrm{G})$ with $E_\mathrm{G}\subseteq E_\mathrm{H}$ and  $\de^\mathrm{G}(a,b)=\de^\mathrm{H}(a,b)$ for all edges $\{a,b\}\in E_\mathrm{G}$.

If $\mathrm{G}$ is connected, regular and every cycle of\/  $\mathrm{G}$ is metric then $\mathrm{M}=(G,\boldsymbol{d}^\mathrm{G})$ is a metric space on $G$ with $\de^\mathrm{G}(a,b)=\boldsymbol{d}^\mathrm{G}(a,b)$ for all edges $\{a,b\}\in E_\mathrm{G}$. 
\end{thm}
\begin{proof}
If there exists a metric space $\mathrm{M}=(G;\de^\mathrm{M})$ with $\dist(\mathrm{M})\subseteq \mathcal{R}$ and $\mathrm{G}\curlyeqprec\mathrm{M}$ then $\mathrm{G}\curlyeqprec (G; [G]^2, \de^\mathrm{M})$ and the $\mathcal{R}$-graph  $(G; [G]^2, \de^\mathrm{M})$, being metric and regular  according to Corollary \ref{cor;mspmetric},  implies that $\mathrm{G}$ is metric. 

The second assertion follows from Lemma \ref{lem:connext} and the third from Lemma \ref{lem:graphext} together with Lemma \ref{lem:cyclemetric}. 
\end{proof}

\section{A construction}\label{section:construction}

\noindent
For  this section, let:\\
$\mathcal{R}\subseteq \Re_{\geq 0}$ be a closed subset of the non negative reals satisfying  the 4-values condition with $0<r\in \mathcal{R}$.   Let $\mathrm{U}=(U;\de^\mathrm{U})$ and $\mathrm{V}=(V;\de^\mathrm{V})$ be two disjoint, countable metric spaces with $\dist(\mathrm{U})\cup \dist(\mathrm{V})\subseteq \mathcal{R}$. Let  $(u_i; i\in \omega)$ be an enumeration of $U$. Let $I\subseteq \omega$ and  $V=\{v_{i}\mid i\in I\}$  an indexing of  $V$  for which $|\de^\mathrm{U}(u_i,u_j)-\de^\mathrm{W}(v_i,v_j)|\leq r$ for all $i,j\in \omega$.

\begin{lem}\label{lem:4cycl}
There exists a metric space $\mathrm{W}=(W;\de^{\mathrm{W}})$ with $\dist(\mathrm{W})\subseteq \mathcal{R}$ and $V\subseteq W=\{w_i;i\in \omega\}$ and $w_i =v_i$ for all $i\in I$  and $W\cap U=\emptyset$, so that $|\de^\mathrm{U}(u_i,u_j)-\de^\mathrm{W}(w_i,w_j)|\leq r$ for all $i,j\in \omega$. 

\end{lem}
\begin{proof}
Let $\mathrm{G}=(G;E_\mathrm{G}, \de^\mathrm{G})$ be the $\mathcal{R}$-graph determined by:
\begin{enumerate}
\item $G=U\cup V$.
\item  $\restrict{\mathrm{G}}{U}=\mathrm{U}$ and $\restrict{\mathrm{G}}{V}=\mathrm{{V}}$. 
\item $E_\mathrm{G}=[U]^2\cup [V]^2\cup \big\{\{u_i,v_i\}\mid i\in I \big\}$.
\item $\de^\mathrm{G}(u_i,v_i)=r$ for all $i\in I$. 
\end{enumerate}
\vskip 5pt

\noindent
Claim: Every cycle of $\mathrm{G}$ is metric. Let $C=\{c_i\mid i\in n\in \omega\}$ induce a cycle of $\mathrm{G}$.  If $C\subseteq U$ or $C\subseteq V$ then, because the only cycles in metric spaces are triangles, the set  $C$ induces a triangle which is metric. If $C\cap U\not=\emptyset$ and $C\cap V\not=\emptyset$ then $|C\cap U|=|C\cap V|=2$. (Every subset $P$ with $|P|<n$ of a cycle of length $n$ induces a set of paths and hence a subgraph with fewer than $|P|$ edges. Every finite subset of $U$ or $V$ induces as many edges as it has elements.) Hence if $C\cap U\not=\emptyset$ and $C\cap V\not=\emptyset$ then there are indices $i,j\in I$ with $i\not=j$ and $C=\{u_{i},u_{j}, v_{j},v_{i}\}$. This cycle is metric because $|\de^\mathrm{U}(u_{i},u_{j})-\de^\mathrm{V}(v_{i},v_{j}))|\leq r$. 

The $\mathcal{R}$-graph $\mathrm{G}$ is regular because if $a$ and $b$ in $G$ are not adjacent then one of them is in $V$ and the other in $U$ and hence every walk from $a$ to $b$ contains an edge of length $r$. It follows from Theorem \ref{thm:graphext} that there exists a metric space $\mathrm{M}=(M=U\cup V;\de)$ with $\mathrm{G}\curlyeqprec\mathrm{M}$ and $\dist(\mathrm{M})\subseteq \mathcal{R}$. Let $W=\{w_i\mid i\in \omega\}$ be a set of points with $w_i=v_i$ for all $i\in I$ and $W\cap U=\emptyset$. Let  $f: W\to M$ be  given by:
\[
f(w_i)=
\begin{cases}
v_i      & \text{for all $i\in I$}, \\
 u_i     & \text{otherwise}.
\end{cases}
\]
Let $\mathrm{W}=(W;\de^\mathrm{W})$ be the metric space on $W$ for which $f$ is an isometry of\/ $\mathrm{W}$ into $\mathrm{M}$. Let $i\in I$ and $j\in \omega\setminus I$ then $|\de^\mathrm{U}(u_i,u_j)-\de^\mathrm{W}(w_i,w_j)|=|\de^\mathrm{M}(u_i,u_j)-\de^\mathrm{M}(v_i,u_j)|\leq r$ because  $(u_i,v_i,u_j)$ is a metric triangle in the metric space $\mathrm{M}$. Hence $|\de^\mathrm{U}(u_i,u_j)-\de^\mathrm{W}(w_i,w_j)|\leq r$ for all $i,j\in \omega$.
\end{proof}

\noindent
We fix the metric space $\mathrm{W}$  given by Lemma \ref{lem:4cycl}.

\vskip 4pt
For $n\in \omega$  let $\boldsymbol{P}_n$ be the set of all order preserving injections $\alpha: \{0,1,2,\dots,n\}\to \omega$ for which the function $\underline{\alpha}: \{u_i\mid 0\leq i\leq n\}\to U$ with $\underline{\alpha}(u_i)=u_{\alpha(i)}$ is an isometry. That is $\de(u_{\alpha(i)},u_{\alpha(j)})=\de(u_i,u_j)$ for all $0\leq i,j\leq n$. Let $\boldsymbol{P}=\bigcup_{n\in\omega} \boldsymbol{P}_n$.    Let $\boldsymbol{P}_\omega$ be the set of all order preserving injections of $\alpha:\omega \to\omega$ for which the function $\underline{\alpha}: U\to U$ with $\underline{\alpha}(u_i)=u_{\alpha(i)}$ is an isometry. Note that $\boldsymbol{P}_\omega$ is the set of order preserving functions $\alpha$ from $\omega$ into $\omega$ for which the set $\{u_{\alpha(i)}=\underline{\alpha}(u_i)\mid i\in \omega\}$ induces a copy of\/ $\mathrm{U}$ in $\mathrm{U}$. This isometric copy of\/  $\mathrm{U}$ in $\mathrm{U}$ is denoted by $\mathrm{U}_\alpha=(U_\alpha;\de^\mathrm{U})$.

For $\alpha\in \boldsymbol{P}_n$  and $\beta\in \boldsymbol{P}_m$ with $ n\leq m\in \omega$ let  $\alpha\sqsubseteq \beta$ if $\alpha(i)=\beta(i)$ for all $i\in n$.  It follows that the partial order $(\boldsymbol{P}; \sqsubseteq )$ is a disjoint union of  trees each one of which having one element of $\boldsymbol{P}_0$ as a root.   If $\boldsymbol{B}\subseteq \boldsymbol{P}$ forms a maximal branch of $(\boldsymbol{P}; \sqsubseteq )$ then $\bigcup\boldsymbol{B}$ is a function in $\boldsymbol{P}_\omega$. If $\alpha\in \boldsymbol{P}_\omega$ let $\alpha_n\in \boldsymbol{P}_n$ be the element with $\alpha_n\sqsubseteq \alpha$.  Then $\boldsymbol{B}=\{\alpha_n\mid n\in \omega\}$ is a maximal branch.  

Let $\mathrm{P}=(\boldsymbol{P};E_\mathrm{P},\de^\mathrm{P})$ be the $\mathcal{R}$-graph with:
\begin{enumerate}
\item $E_\mathrm{P}=\big\{(\alpha,\beta)\in [\boldsymbol{P}]^2\mid \alpha\sqsubseteq \beta\big\}$. 
\item $\de^\mathrm{P}(\alpha,\beta)=\de^\mathrm{W}(w_{n},w_{m})$  for all  $\alpha\in \boldsymbol{P}_n$ and $\beta\in \boldsymbol{P}_m$ with $\alpha\sqsubseteq \beta$.
\end{enumerate}  
Note that if $\boldsymbol{B}$ is a  branch of $(\boldsymbol{P};\sqsubseteq)$ then $\boldsymbol{B}$ induces a metric space in $\mathrm{P}$ and if $\boldsymbol{B}$ is an infinite maximal branch then $\boldsymbol{B}$ induces an isometric copy of $\mathrm{W}$ in $\mathrm{P}$. (It is not difficult to see that every maximal branch is infinite, but we do not need this fact.)

\begin{lem}\label{lem:Pmetric}
Every induced cycle of the $\mathcal{R}$-graph $\mathrm{P}=(\boldsymbol{P};E_\mathrm{P},\de^\mathrm{P})$ is metric.
\end{lem}
\begin{proof}
Let $C$ induce a cycle in $\mathrm{P}$ and $n$ maximal so that there is a point $\alpha\in C\cap \boldsymbol{P}_n$. Let $\beta$ and $\gamma$ be the two points in $C$ adjacent to $\alpha$. It follows that $n\geq 2$ and that  $\beta$ and $\gamma$ are on the branch of the tree $(\boldsymbol{P};\sqsubseteq)$ below $\alpha$ and hence adjacent. Hence $C=\{\alpha,\beta,\gamma\}$ is a triangle in the metric space induced by the branch below $\alpha$. 
\end{proof}

Let $\mathrm{H}=(H;E_\mathrm{H},\de^\mathrm{H})$ be the $\mathcal{R}$-graph with $H=\boldsymbol{P}\cup U$ and $E_\mathrm{H}=E_\mathrm{P}\cup [U]^2\cup \bigcup_{n\in \omega}\big\{\{\alpha,u_{\alpha(n)}\}\mid \alpha\in \boldsymbol{P}_n\big\}$ so that:
\begin{enumerate}
\item $\restrict{\mathrm{H}}{\boldsymbol{P}}=\mathrm{P}$ and $\restrict{\mathrm{H}}{U}=\mathrm{U}$. 
\item $\de^\mathrm{H}(\alpha,u_{\alpha(n)})=r$ for all $n\in \omega$ and $\alpha\in \boldsymbol{P}_n$. 
\end{enumerate}
Note that in $\mathrm{H}$ every point $\alpha\in \boldsymbol{P}$ is adjacent to exactly one point in $U$. Namely, if $\alpha\in \boldsymbol{P}_n$ then it is adjacent to $u_{\alpha(n)}$.  Let $u^{\langle\alpha\rangle}$ be  the point in $U$ adjacent to $\alpha\in \boldsymbol{P}$. Then $\de^\mathrm{H}(\alpha,u^{\langle\alpha\rangle})=r$ and for $\alpha \sqsubseteq \beta$ two points in $\boldsymbol{P}$: 
\begin{align}\label{align:redin}
|\de^\mathrm{H}(\alpha,\beta)-\de^\mathrm{H}(u^{\langle\alpha\rangle},u^{\langle\beta\rangle})|\leq r, 
\end{align}
because for $\alpha\in \boldsymbol{P}_n$ and $\beta\in \boldsymbol{P}_m$:  \\ 

\vskip -6pt
\noindent
$|\de^\mathrm{H}(\alpha,\beta)-\de^\mathrm{H}(u^{\langle\alpha\rangle},u^{\langle\beta\rangle})|=
|\de^\mathrm{W}(w_{n},w_{m})-\de^\mathrm{U}(u_{\alpha(n)},u_{\alpha(m)})|= \\
|\de^\mathrm{W}(w_{n},w_{m})-\de^\mathrm{U}(u_{n},u_{m})|\leq r$. Inequality \ref{align:redin} implies that:
\begin{align}\label{align:redin1}
\de^\mathrm{H}(\alpha,\beta)\leq r\oplus \de^\mathrm{H}(u^{\langle\alpha\rangle},u^{\langle\beta\rangle}) \text{  and $\de^\mathrm{H}(u^{\langle\alpha\rangle},u^{\langle\beta\rangle})\leq r\oplus \de^\mathrm{H}(\alpha,\beta)$}.
\end{align}

\begin{lem}\label{lem:bupbl}
Let $\mathrm{U}^\ast=(U^\ast;\de^\mathrm{U})$ be an isometric copy of\/ $\mathrm{U}$ in $\mathrm{U}$. Then there exists an isometric copy $\mathrm{V}^\ast=(\mathrm{V}^\ast;\de)$ of\/ $\mathrm{V}$ in $\mathrm{H}$ with:
\[
V^\ast\subseteq \big(U^\ast\big)_r.
\]
\end{lem}
\begin{proof}
It follows from Lemma \ref{lem:ordemb} that there exists an enumeration order preserving isometry $\underline{\alpha}$ of\/ $\mathrm{U}$ into $\mathrm{U}$ with $\alpha(U)\subseteq U^\ast$. Let $\alpha\in \boldsymbol{P}_\omega$ be the order preserving function of $\omega$ into $\omega$ with $\underline{\alpha}(u_i)=u_{\alpha(i)}$ and let $\alpha_n$ be the restriction of $\alpha$ to $n$. then $\boldsymbol{B}=\{\alpha_n\mid n\in \omega\}$ is an infinite maximal chain in $(\boldsymbol{P};\sqsubseteq)$ and hence isometric to $\mathrm{W}$. Because $\mathrm{V}$ is a subspace of $\mathrm{W}$ there exists an isometric copy $\mathrm{V}^\ast=(V^\ast;\de)$ of \/ $\mathrm{V}$ in $\boldsymbol{B}$. For every point $\beta\in \boldsymbol{B}$ the distance $\de^\mathrm{H}(\beta,u^{\langle\beta\rangle})=r$ and hence: $V^\ast\subseteq \boldsymbol{B}\subseteq \big(\alpha(U)\big)_r\subseteq \big(U^\ast\big)_r$.
\end{proof}

\begin{lem}\label{lem:Hmetric}
The $\mathcal{R}$-graph $\mathrm{H}$ is metric and regular.
\end{lem}
\begin{proof}
Let $C$ induce a cycle in $\mathrm{H}$. It follows from Lemma \ref{lem:Pmetric} that if $C\subseteq \boldsymbol{P}$ then $C$ is metric. If $L$ induces a metric subspace of $\mathrm{H}$ and $C\subseteq L$ then $C$ is a metric  triangle. Hence we may assume that $C$ is not a subset of any metric subspace of $\mathrm{H}$. It follows that $|C\cap L|\leq 2$ for every subset $L$ of $H$ which induces a metric subspace of $\mathrm{H}$. In particular $1\leq |C\cap U|\leq 2$ and $|C\cap \boldsymbol{B}|\leq 2$ for every branch $\boldsymbol{B}$ of $\mathrm{P}$.  

Let $\{\alpha_i\mid i\in l\in \omega\}$ be the set of points in $C\cap \boldsymbol{P}$ for which there does not exist a point $\beta\in C\cap \boldsymbol{P}$ with $\alpha_i \sqsubset \beta$, that is those points in $C\cap \boldsymbol{P}$ which are maximal on their branches. Let $\boldsymbol{B}_i$ the set of points $\beta\in \boldsymbol{P}$ with $\beta\sqsubset \alpha_i$. Then $|C\cap \boldsymbol{B}_i|\leq 2$. Hence for each $\alpha_i$ there is at most on other point, say $\beta_i$ in $C\cap \boldsymbol{P}$, which implies that the other point of $C$ adjacent to $\alpha_i$ is the point $u^{\langle\alpha_i\rangle}$ in $U$. That is the two points in $C$ adjacent to $\alpha_i$ are  $u^{\langle\alpha_i\rangle}$ and $\beta_i$. 

Let $c\in C$ be the other point of $C$ adjacent to $\beta_0$. If $c\in \boldsymbol{P}$ then it is not possible that $c\sqsubset \beta_0$ because otherwise $\alpha_0$ is adjacent to $c$. Hence if $c\in \boldsymbol{P}$ then $\beta_0 \sqsubset c$.  If there is an $\alpha_i$ with $c\sqsubset \alpha_i$ the branch below or equal to $\alpha_i$ would contain three points of $C$. Hence $c=\alpha_i$ for some $i\in l$. Say $i=1$ and then the other point of $C$ adjacent to $\beta_0$ is $\alpha_1$. That is, $\beta_0=\beta_1$. Similar to before, the other point of $C$ adjacent to $\alpha_1$ is the point $u^{\langle\alpha_1\rangle}$ in $U$. The circle $C$ contains the path $u^{\langle\alpha_1\rangle}, \alpha_1,\beta_0,\alpha_0,u^{\langle\alpha_0\rangle}$. If $u^{\langle\alpha_1\rangle}\not= u^{\langle\alpha_0\rangle}$ then $C=\{u^{\langle\alpha_1\rangle}, \alpha_1,\beta_0,\alpha_0,u^{\langle\alpha_0\rangle}\}$  because $\{u^{\langle\alpha_1\rangle}, u^{\langle\alpha_0\rangle}\}$ is an edge in $\mathrm{H}$. In order to see that $C$ induces a metric cycle we have to check for each edge $\{x,y\}$ induced by $C$ that $\de(x,y)$ is less than or equal to the $\bigoplus$-sum over the other edges. Because $\de(\alpha_0, u^{\langle\alpha_0\rangle})=\de(\alpha_1, u^{\langle\alpha_1\rangle})$ this will be the case for the edges $\{(\alpha_0, u^{\langle\alpha_0\rangle}\}$ and $\{(\alpha_0, u^{\langle\alpha_0\rangle}\}$. Using Inequality \ref{align:redin1}:
\begin{align*}
&\de^\mathrm{H}(u^{\langle\alpha_0\rangle},u^{\langle\alpha_1\rangle})\leq \de^\mathrm{H}(u^{\langle\alpha_0\rangle},u^{\langle\beta_0\rangle})\oplus \de^\mathrm{H}(u^{\langle\beta_0\rangle},u^{\langle\alpha_1\rangle})\leq \\
&\big(r\oplus \de^\mathrm{H}(\alpha_0,\beta_0)\big)\oplus\big(r\oplus\de^\mathrm{H}(\beta_0,\alpha_1)\big)
\end{align*}
and
\begin{align*}
&\de^\mathrm{H}(\alpha_0,\beta_0)\leq r\oplus \de^\mathrm{H}(u^{\langle\alpha_0\rangle},u^{\langle\beta_0\rangle})\leq r\oplus \de^\mathrm{H}(u^{\langle\beta_0\rangle},u^{\langle\alpha_1\rangle})\oplus \de^\mathrm{H}(u^{\langle\alpha_0\rangle},u^{\langle\alpha_1\rangle})\\
&\leq r\oplus \de^\mathrm{H}(\beta_0,\alpha_1) \oplus r\oplus \de^\mathrm{H}(u^{\langle\alpha_0\rangle},u^{\langle\alpha_1\rangle})
\end{align*}
and similar for $\de^\mathrm{H}(\alpha_1,\beta_0)$.

The remaining possibility is that $c\in U$. Then $C$ contains the path $u^{\langle\alpha_0\rangle}, \alpha_0,\beta_0, u^{\langle\beta_0\rangle}$ and hence $C=\{u^{\langle\alpha_0\rangle}, \alpha_0,\beta_0, u^{\langle\beta_0\rangle}\}$  because $u^{\langle\beta_0\rangle}$ and $u^{\langle\alpha_0\rangle}$ are adjacent in $\mathrm{H}$. The cycle induced by $C$ is metric because: 
\begin{align*}
&\de^\mathrm{H}(u^{\langle\alpha_0\rangle}, u^{\langle\beta_0\rangle})\leq \de^\mathrm{H}(\alpha_0,\beta_0)\oplus r\leq\de^\mathrm{H}(\alpha_0,\beta_0)\oplus r\oplus r \text{  and}\\
&\de^\mathrm{H}(\alpha_0,\beta_0)\leq \de^\mathrm{H}(u^{\langle\alpha_0\rangle}, u^{\langle\beta_0\rangle})\oplus r\leq \de^\mathrm{H}(u^{\langle\alpha_0\rangle}, u^{\langle\beta_0\rangle})\oplus r\oplus r.
\end{align*}

Let $a$ and $b$ be two non adjacent points in $H$. If one is an element of $P$ and the other of $U$ then every walk from $a$ to $b$ contains an edge of length $r$. If both are in $U$ then every walk from $a$ to $b$ has length at least $\de^\mathrm{H}(a,b)$. If both points are in $P$ then every walk from $a$ to $b$ which contains a point of $U$ has length at least $r$. The remaining case then is that $a$ and $b$ are two non adjacent points in the same connected component of $\mathrm{P}$. This connected component is a tree. Let $n\in \omega$ maximal with $\gamma\in \boldsymbol{P}_n$ and $\gamma\sqsubset a$ and $\gamma\sqsubset b$. Then $\boldsymbol{d}^\mathrm{H}(a,b)\geq \de^\mathrm{H}(a,\gamma)$.

\end{proof}

\begin{thm}\label{thm:treec}
$\mathcal{R}\subseteq \Re_{\geq 0}$ be a closed subset of the non negative reals satisfying  the 4-values condition with $0<r\in \mathcal{R}$.   Let $\mathrm{U}=(U;\de^\mathrm{U})$ and $\mathrm{V}=(V;\de^\mathrm{V})$ be two disjoint metric spaces with $\dist(\mathrm{U})\cup \dist(\mathrm{V})\subseteq \mathcal{R}$ and let  $(u_i; i\in \omega)$ be an enumeration of $U$. Let $I\subseteq \omega$ and  $V=\{v_{i}\mid i\in I\}$ be an indexing of  $V$.  Let 
\[
|\de^\mathrm{U}(u_i,u_j)-\de^\mathrm{V}(v_i,v_j)|\leq r  \text{  \  for all $i,j\in I$. }
\]
Then there exists a countable metric space $\mathrm{L}=(L;\de^\mathrm{L})$ with $\dist(\mathrm{L})\subseteq \mathcal{R}$, containing the metric space $\mathrm{U}$ as a subspace, so that for every copy $\mathrm{U}^\ast=(U^\ast;\de^\mathrm{L})$ of\/ $\mathrm{U}$ in $\mathrm{U}$ there exists a copy $\mathrm{V}^\ast=(V^\ast;\de^\mathrm{L})$ of\/ $\mathrm{V}$ in $\mathrm{L}$ so that:
\[
V^\ast\subseteq \big(U^\ast\big)_r.
\]
\end{thm}
\begin{proof}
The $\mathcal{R}$-graph $\mathrm{H}$ is metric according to  Lemma \ref{lem:Hmetric} and hence it follows from Theorem \ref{thm:graphext} that there is a metric space $\mathrm{L}=(L;\de)$ with   $\mathrm{H}\curlyeqprec \mathrm{L}$ and $\dist(\mathrm{L})\subseteq \mathcal{R}$ and $L=U\cup \boldsymbol{P}$. Lemma \ref{lem:bupbl} implies that for every copy $\mathrm{U}^\ast=(U^\ast;\de^\mathrm{L})$ of\/ $\mathrm{U}$ in $\mathrm{U}$ there exists a copy $\mathrm{V}^\ast=(V^\ast;\de^\mathrm{L})$ of\/ $\mathrm{V}$ in $\mathrm{L}$ so that:
\[
V^\ast\subseteq \big(U^\ast\big)_r.
\]

\end{proof}

\section{Finite approximations}\label{sect:finappr}

Let $\mathcal{R}\subseteq \Re_{\geq 0}$ be a closed and bounded subset of the non negative reals satisfying  the 4-values condition. We will write $\oplus$ for $\oplus_\mathcal{R}$.  Let $\epsilon>0$ be given.

\begin{defin}\label{defin:rooted}
For $A$ a finite subset of\/ $\mathcal{R}$ with $\max{A}=\max\mathcal{R}$ and $l\in \mathcal{R}$ let $\bar{l}^{\langle A\rangle}$ or just $\bar{l}$, if  $A$ is understood, be the smallest  number in $A$ larger than or equal to $l$. 

A finite set $A\subseteq \mathcal{R}$ is a {\em finite  $\epsilon$-approximation} of\/ $\mathcal{R}$ if  $\max A= \max\mathcal{R}$ and $\bar{l}-l<\epsilon$ for all $l\in \mathcal{R}$ with $l\not=0$. 
\end{defin}
\noindent
Because $\mathcal{R}$ is compact there exists for every finite subset $B\subseteq \mathcal{R}$ and every $\epsilon>0$  a finite $\epsilon$-approximation $A$ of $\mathcal{R}$  with $B\subseteq A$.

\begin{defin}\label{defin:subadd}
For $A\subseteq \mathcal{R}$  let $\mathfrak{C}(A)=\{a\oplus b\mid a,b\in A\} \cup A$ and $\mathfrak{C}(\emptyset)=A$. The set $A$  is {\em subadditive closed } if   $\mathfrak{C}(A)=A$.

Let   $\mathfrak{C}_0(A)=\emptyset $ and $\mathfrak{C}_1(A)=A$ and recursively $\mathfrak{C}_{n+1}(A)=\mathfrak{C}(\mathfrak{C}_n(A))$. 
\end{defin}

Note  that $\min A=\min \mathfrak{C}(A)$ and if $\max A=\max\mathcal{R}$ then  $\max A=\max \mathfrak{C}(A)$ and if $A$ is finite then $\mathfrak{C}(A)$ is finite.  Hence if  $A$ is a finite  $\epsilon$-approximation of $\mathcal{R}$ then $\mathfrak{C}_n(A)$ is a finite $\epsilon$-approximation of $\mathcal{R}$.  If $\mathfrak{C}_n(A)$ is subadditive closed then $\mathfrak{C}_{n+1}(A)=\mathfrak{C}_n(A)$ is subadditive closed.  

For $\mathfrak{C}_n(A)$ not subadditive closed   let $\boldsymbol{w}_{n+1}:=\min\big(\mathfrak{C}_{n+1}(A)\setminus\mathfrak{C}_n(A)\big)$. (Hence $\boldsymbol{w}_1=\min A$.)      

\begin{lem}\label{lem:b1}
Let  $A$ be a finite  $\epsilon$-approximation of $\mathcal{R}$.  Let  $1\leq n\in \omega$ so that $\mathfrak{C}_n(A)$ is not subadditive closed, then:
\begin{enumerate}
\item $\boldsymbol{w}_{n+1}\geq \boldsymbol{w}_n\oplus \boldsymbol{w}_1$.
\item $\boldsymbol{w}_{n+1}>\boldsymbol{w}_n$
\item If $\boldsymbol{w}_{n+1}>\boldsymbol{w}_n\oplus \boldsymbol{w}_1$ then $\boldsymbol{w}_{n+1}\geq \boldsymbol{w}_n+\boldsymbol{w}_1$.  
\end{enumerate} 
\end{lem}
\begin{proof}
There exist $t'\in \mathfrak{C}_{n}(A)\setminus\mathfrak{C}_{n-1}(A)$ and $t''\in \mathfrak{C}_{n}(A)$ so that $\boldsymbol{w}_{n+1}=t'\oplus t''$. Hence $\boldsymbol{w}_{n+1}\geq t'\oplus t''\geq \boldsymbol{w}_n\oplus\boldsymbol{w}_1\geq \boldsymbol{w}_n$. It is not the case that  $\boldsymbol{w}_{n+1}\not=\boldsymbol{w}_n$ because $\boldsymbol{w}_{n+1}\in \mathfrak{C}_{n+1}(A)\setminus \mathfrak{C}_n(A)$ and $\boldsymbol{w}_n\in \mathfrak{C}_n(A)$.  Hence $\boldsymbol{w}_{n+1}>\boldsymbol{w}_n$.

There is no element $x\in \mathcal{R}$ with $\boldsymbol{w}_n\oplus \boldsymbol{w}_1<x\leq\boldsymbol{w}_n+\boldsymbol{w}_1$. Hence if  $\boldsymbol{w}_{n+1}>\boldsymbol{w}_n\oplus \boldsymbol{w}_1$ then $\boldsymbol{w}_{n+1}\geq \boldsymbol{w}_n+\boldsymbol{w}_1$.
\end{proof}

\begin{lem}\label{lem:1}
Let  $A$ be a finite  $\epsilon$-approximation of\   $\mathcal{R}$. There is a number $n\in \omega$ so that $\mathfrak{C}_n(A)$ is a  subadditive closed and finite $\epsilon$-approximation of $\mathcal{R}$.
 \end{lem}
 \begin{proof}
The set $\mathfrak{C}_n(A)$ is an $\epsilon$-approximation of $\mathcal{R}$ for every $1\leq n\in \omega$. Hence,  because $\mathcal{R}$ is bounded above, it suffices to show:  If  $n\geq 1$ and  $\mathfrak{C}_{n+1}(A)$ is not subadditive closed  then 
$\boldsymbol{w}_{n+2} \geq \boldsymbol{w}_n+\boldsymbol{w}_1$.

If $\boldsymbol{w}_n\oplus \boldsymbol{w}_1<\boldsymbol{w}_{n+1}$ then $\boldsymbol{w}_{n+2}\geq \boldsymbol{w}_{n+1}\geq \boldsymbol{w}_n+\boldsymbol{w}_1$ according to Lemma \ref{lem:b1}. Because $\boldsymbol{w}_n\oplus \boldsymbol{w}_1\leq \boldsymbol{w}_{n+1}$ the remaining case is $\boldsymbol{w}_n\oplus \boldsymbol{w}_1=\boldsymbol{w}_{n+1}$.

Then using Lemma \ref{lem:b1}: $\boldsymbol{w}_{n+2}>\boldsymbol{w}_{n+1}=\boldsymbol{w}_n\oplus \boldsymbol{w}_1$. There is no element $x\in \mathcal{R}$ with $\boldsymbol{w}_n\oplus \boldsymbol{w}_1<x< \boldsymbol{w}_n+\boldsymbol{w}_1$ implying $\boldsymbol{w}_{n+2}\geq \boldsymbol{w}_n+\boldsymbol{w}_1$. 
 \end{proof}

 \begin{cor}\label{cor:111}
 There exists for every $\epsilon>0$   and finite $B\subseteq \mathcal{R}$    a finite $\epsilon$-approximation $\mathcal{S}$  of $\mathcal{R}$ which is subadditive closed and with $B\subseteq \mathcal{S}$. 
 
If $0$ is a limit of $\mathcal{R}$, then there exist, for every $\epsilon>0$ and $\mathcal{R}\ni r<\epsilon$   and finite $B\subseteq \mathcal{R}$    a finite $r$-approximation $\mathcal{S}$  of $\mathcal{R}$ which is subadditive closed and with $B\subseteq \mathcal{S}$ and with $r$ the minimum of $\mathcal{S}$.  
\end{cor}

\begin{lem}\label{lem:add}
Let $A\subseteq \mathcal{R}$ be finite and  subadditive closed with $\max A=\max\mathcal{R}$ and for $x\in A$ let $\bar{x}:=\bar{x}^{\langle A\rangle}$.   Then if   $\{a,b,c\}\subseteq \mathcal{R}$ is metric the set $\{\bar{a},\bar{b},\bar{c}\}$ is metric.
\end{lem}
\begin{proof}
Let $\{a,b,c\}$ be metric with  $a\leq b\leq c$. Then $\bar{a}+\bar{b}\geq a+b\geq c$ and $c\in \mathcal{R}$  implies $\bar{a}\oplus \bar{b}\geq c$. This in turn implies that  $\bar{c}\leq \bar{a}\oplus \bar{b}$ because $\bar{a}\oplus \bar{b} \in \mathcal{R}$. Finally  $\bar{c}\leq \bar{a}\oplus \bar{b}\leq \bar{a}+\bar{b}$. 
 
\end{proof}

\begin{lem}\label{lem:add4val}
Let $A\subseteq \mathcal{R}$ be finite and  subadditive closed with $\max A=\max\mathcal{R}$. Then $A$ satisfies the 4-values condition.
\end{lem}
\begin{proof}
For $l\in \mathcal{R}$ let $\bar{l}:=\bar{l      }^{\langle A\rangle}$.   
 
 Let $(a,b,c,d)\in \mathscr{Q}(A)$ with $x\in A$ so that $x\leadsto (a,b,c,d)$. Because $A\subseteq \mathcal{R}$, and $A$ satisfies the 4-values condition,  there is a number $y\in \mathcal{R}$ with $y\leadsto (a,d,c,b)$ implying that the triples $(y,c,b)$ and $(y,a,d)$ are metric. It follows from Lemma~\ref{lem:add} that then the triples $(\bar{y},\bar{c}, \bar{b})=(\bar{y},c,b)$ and $(\bar{y},\bar{a},\bar{d})=(\bar{y},a,d)$ are metric and hence that $\bar{y}\leadsto (a,d,c,b)$.  
 
\end{proof}

\begin{thm}\label{thm:approxfin}
Let $\mathcal{R}$ be a closed and bounded subset of $\Re_{\geq 0}$ satisfying the 4-values condition. Let  $\UR=(U_R;\de^{\UR})$ be the Urysohn metric space with with $\dist(\UR)=\mathcal{R}$ and let $V$ be a countable  subset of $\UR$.  

Then there exists, for every $\epsilon>0$,  a finite subset $\mathcal{S}$ of $\mathcal{R}$, which satisfies the 4-values condition, and there exists an injection $f$ of  $V$ into $U$, with $\boldsymbol{U}=(U;\de^{\boldsymbol{U}})$ a copy of the Urysohn metric space with $\dist(\boldsymbol{U})=\mathcal{S}$, so that for all $a,b\in V$:
\begin{align}\label{align:diff17}
0<\de^{\boldsymbol{U}}(f(a),f(b))-  \de^{\UR}(a,b)<\epsilon.
\end{align}
For every finite $B\subseteq \mathcal{R}$ the set $\mathcal{S}$ can be constructed to contain $B$ as a subset. The Urysohn metric space $\boldsymbol{U}$ can be chosen with $U\cap U_\mathcal{R}=\emptyset$. If 0 is a limit of $\mathcal{R}$  and for $\mathcal{R}\ni r<\epsilon$ the set $\mathcal{S}$ can be taken to contain $r$ as a minimum and the number $\epsilon$ in Inequality \ref{align:diff17} can be replaced by $r$. 
\end{thm}
\begin{proof}
There exists, according to Corollary \ref{cor:111}, a finite $\epsilon$-approximation $\mathcal{S}$  of $\mathcal{R}$ which is subadditive closed and with $B\subseteq \mathcal{S}$. It follows from Lemma \ref{lem:add4val} that $\mathcal{S}$ satisfies the 4-values condition. Let $\boldsymbol{U}=(U;\de^{\boldsymbol{U}})$ be the Urysohn metric space with $\dist(\boldsymbol{U})=\mathcal{S}$ given by Theorem \ref{thm:characterization}. For every pair of points $p,q\in U_\mathcal{R}$ and with $l=\de^{\UR}(p,q)$ let $\bar{\de}(p,q)=\bar{l}^{\langle\mathcal{S}\rangle}$. It follows from Lemma \ref{lem:add} that $\bar{\de}$ is a metric on $U_\mathcal{R}$ and  
\begin{align}\label{align:near}
0<\bar{\de}(p,q)-\de^{\UR}(p,q)<\epsilon,
\end{align}
because $\mathcal{S}$ is an $\epsilon$-approximation of $\mathcal{R}$. Let $\mathrm{V}$ be the metric space with $V$ as set of points and $\bar{\de}$ as metric function. It follows from Fact \ref{fact:first} that there is an isometric embedding $f$ of the metric space $\mathrm{V}$ into the Urysohn space $\boldsymbol{U}$. From Inequality \ref{align:near}:
\[
0<\de^{\boldsymbol{U}}(f(a),f(b))-  \de^{\UR}(a,b)<\epsilon.
\]
\end{proof}

\section{Oscillation stability theorems}\label{sect:oscillstab}

\begin{lem}\label{thm:fext}
Let $\mathcal{R}\subseteq \Re_{\geq 0}$ be a closed and bounded subset of the non negative reals satisfying  the 4-values condition and with 0 as a limit.  Let $\UR=(U_\mathcal{R};\de)$ be an Urysohn metric space with $\dist(\UR)=\mathcal{R}$ and let $V$ be a countable dense subset of $U_\mathcal{R}$ and $\restrict{\UR}{V}$ be equal to the space $\mathrm{V}=(V;\de)$. 

Then: There exists,  for every $\epsilon>0$,  a finite set $\mathcal{S}\subseteq \mathcal{R}$ satisfying the 4-values condition and a subset $U_\mathcal{S}$ of $U_\mathcal{R}$ so that $\restrict{\UR}{U_\mathcal{S}}$ is the Urysohn metric space $\Ur{\mathcal{S}}$. For every isometric copy $\Ur{\mathcal{S}}^\ast=(U^\ast_\mathcal{S}, \de)$ of $\Ur{\mathcal{S}}$ in $\Ur{\mathcal{S}}$ there exists an isometric copy $\mathrm{V}_{\Ur{\mathcal{S}}^\ast}=(V_{\Ur{\mathcal{S}}^\ast};\de)$ of\/ $\mathrm{V}$ in $\UR$ with:
\[
V_{\Ur{\mathcal{S}}^\ast}\subseteq \big(U^\ast_\mathcal{S}\big)_\epsilon.
\]
\end{lem}
\begin{proof}
Let $\epsilon>0$ be given and let $0<r\in \mathcal{R}$ with $r<\epsilon$. Theorem~\ref{thm:approxfin} supplies a finite subset $\mathcal{S}\subseteq \mathcal{R}$ with $r$ the minimum of $\mathcal{S}$ and a copy $\boldsymbol{U}=(U;\de^{\boldsymbol{U}})$ of the Uryosohn metric space with $\dist(\boldsymbol{U}))=\mathcal{S}$ and $U\cap U_\mathcal{R}=\emptyset$. Let $f$ be the injection of $V$ into $U$ given by  Theorem~\ref{thm:approxfin} for which 
\begin{align}\label{align:diff171}
0<\de^{\boldsymbol{U}}(f(a),f(b))-  \de^{\UR}(a,b)<r  \text{   \,  for all $a,b\in V$}.
\end{align}

Let $(u_i; i\in \omega)$ be an enumeration of $U$. Let $I=\{i\in \omega\mid \exists v\in V\, (f(v)=u_i)\}$. Let $V=\{v_i\mid i\in I\}$ be the indexing of $V$ with $I$ so that $f(v_i)=u_i$ for all $i\in I$. Let $\de^{\mathrm{V}}$ be the restriction of the metric $\de$ to the subset $V$ of $U_\mathcal{R}$.  It follows then from Inequality \ref{align:diff171} that:
\[
|\de^\mathrm{U}(u_i,u_j)-\de^\mathrm{V}(v_i,v_j)|\leq r  \text{  \  for all $i,j\in I$. }
\] 
We are now in the position to apply  Theorem \ref{thm:treec} which yields a countable metric space $\mathrm{L}=(L;\de^\mathrm{L})$ with $\dist(\mathrm{L})\subseteq \mathcal{R}$, containing the metric space $\mathrm{U}$ as a subspace, so that for every copy $\mathrm{U}^\ast=(U^\ast;\de^\mathrm{L})$ of\/ $\mathrm{U}$ in $\mathrm{U}$ there exists a copy $\mathrm{V}^\ast=(V^\ast;\de^\mathrm{L})$ of\/ $\mathrm{V}$ in $\mathrm{L}$ so that:
\[
V^\ast\subseteq \big(U^\ast\big)_r.
\]
There exists an isometric embedding $g$ of $\mathrm{L}$ into $\UR$, according to Fact~\ref{fact:first}. Let $U_\mathcal{S}$ be the image of $U$ under $g$. Then $\UR$ restricted to  $U_\mathcal{S}$ is an Urysohn space $\Ur{\mathcal{S}}$. If $\Ur{\mathcal{S}}^\ast=(U^\ast_\mathcal{S}, \de)$ is an isometric copy of  of $\Ur{\mathcal{S}}$ in $\Ur{\mathcal{S}}$, the inverse isometry of $g$ yields an isometric copy of\/ $\mathrm{U}$ in $ \mathrm{U}$ and hence a copy of $V$ which projected via $g$ into $U_\mathcal{R}$ gives a copy $\mathrm{V}_{\Ur{\mathcal{S}}^\ast}=(V_{\Ur{\mathcal{S}}^\ast};\de)$  of\/  $\mathrm{V}$ in $\UR$ with:
\[
V_{\Ur{\mathcal{S}}^\ast}\subseteq \big(\Ur{\mathcal{S}}^\ast\big)_r\subseteq \big(U^\ast_\mathcal{S}\big)_\epsilon.
\]

\end{proof}

\begin{thm}\label{thm:final}
Every  bounded, separable, complete, homogeneous, universal metric space $\UR$ with $0$ a limit of $\mathcal{R}$  is oscillation stable.
\end{thm}
\begin{proof}
On account of Theorem \ref{thm:appequosc} it suffices to prove that $\UR=(U_\mathcal{R};\de)$ is approximately indivisible. Let $\epsilon>0$ and the function $\gamma:U_\mathcal{R}\to n\in \omega$ be given and let $V$ be a dense subset of $U_\mathcal{R}$. Lemma \ref{thm:fext} yields a finite set $\mathcal{S}\subseteq \mathcal{R}$ satisfying the 4-values condition  for which  the Urysohn space $\Ur{\mathcal{S}}$ is a subspace of $\UR$ with: For every isometric copy $\Ur{\mathcal{S}}^\ast=(U^\ast_\mathcal{S}, \de)$ of $\Ur{\mathcal{S}}$ in $\Ur{\mathcal{S}}$ there exists an isometric copy $\mathrm{V}_{\Ur{\mathcal{S}}^\ast}=(V_{\Ur{\mathcal{S}}^\ast};\de)$ of\/ $\mathrm{V}$ in $\UR$ with:
\[
V_{\Ur{\mathcal{S}}^\ast}\subseteq \big(U^\ast_\mathcal{S}\big)_{\frac{\epsilon}{2}}.
\]
The restriction of $\gamma$ to $U_\mathcal{S}$ maps $U_\mathcal{S}$ to $n$. Because $\Ur{\mathcal{S}}$ is indivisible according to Theorem \ref{thm:finitedist} there exists $i\in n$ and an isometric copy $\Ur{\mathcal{S}}^\ast=(U^\ast_\mathcal{S};\de)$ of $\Ur{\mathcal{S}}$ in $\Ur{\mathcal{S}}$ with $U^\ast_\mathcal{S}\subseteq \gamma^{-1}(i)$ and hence a copy $\mathrm{V}_{\Ur{\mathcal{S}}^\ast}=(V_{\Ur{\mathcal{S}}^\ast};\de)$ of\/ $\mathrm{V}$ in $\UR$ so that 
\[
V_{\Ur{\mathcal{S}}^\ast}\subseteq \big(U^\ast_\mathcal{S}\big)_{\frac{\epsilon}{2}}\subseteq \big(\gamma^{-1}(i)\big)_{\frac{\epsilon}{2}}.
\]
The completion $C$ of $V_{\Ur{\mathcal{S}}^\ast}$ in $\UR$ induces in $\UR$ a subspace $\mathrm{C}=(C;\de)$ isometric to $\UR$ with:
\[
C\subseteq \big(V_{\Ur{\mathcal{S}}^\ast}\big)_{\frac{\epsilon}{2}}\subseteq \left(\big(\gamma^{-1}(i)\big)_{\frac{\epsilon}{2}}\right)_{\frac{\epsilon}{2}}\subseteq \big(\gamma^{-1}(i)\big)_\epsilon.
\]
\end{proof}
Let $\UR$ be a bounded, uncountable, separable, complete,  homogeneous, universal  metric space.  Because $\UR$ is separable, 0 is a limit of the set of distances $\mathcal{R}$. Hence we obtain:

\begin{thm}\label{thm:finabstr}
Every bounded, uncountable, separable, complete,  homogeneous, universal  metric space is oscillation stable.
\end{thm}

\section{Oscillation stable and approximately indivisible}

It follows from a general discussion of ``oscillation stable" topological groups in \cite{KPT} and \cite{Pe1} that a homogeneous metric space  $\mathrm{M}=(M; \de)$ is  oscillation stable if and only if it is   approximately indivisible. We will prove in this section the more general fact that a  metric space  $\mathrm{M}=(M; \de)$ is  oscillation stable if and only if it is   approximately indivisible.

\begin{lem}\label{claim:1}
If\/  $\mathrm{M}=(M; \de)$ is an   approximately indivisible metric space then $\mathrm{M}$ is oscillation stable.
\end{lem}
\begin{proof}
Let $f: M\to \Re$ be bounded and uniformly continuous and let $\epsilon>0$ be given.   

Let $\delta>0$ be such that $\de(p,q)<\delta$ implies $|f(p)-f(q)|<\frac{1}{3}\epsilon$ for all $p,q\in M$. Let $a,b\in\Re$ be such that $f(M)\subseteq [a,b)$ and let $a=x_0<x_1<x_2<\dots<x_n=b$ be a partition of $[a,b)$ with $x_{i+1}-x_i<\frac{1}{3}\epsilon$ for all $i\in n$. Let   $\gamma: M\to n$ be the function with  $\gamma(p)=i$ if $f(p)\in [x_i,x_{i+1})$.

Because $\mathrm{M}$ is approximately indivisible there exists an $i\in n$ and  a copy $\mathrm{M}^\ast=(M^\ast, \de)$ of $\mathrm{M}$ in $\mathrm{M}$ with $M^\ast\subseteq \big(\gamma^{-1}(i)\big)_\delta$.  Let $p,q\in M^\ast$. There exist $u,v\in \gamma^{-1}(i)$ with $\de(p,u)<\delta$ and $\de(q,v)<\delta$. Then: 
\[
|f(p)-f(q)|\leq |f(p)-f(u)|+|f(u)-f(v)|+|f(v)-f(q)|=\epsilon.
\]
\end{proof}

\begin{lem}\label{claim:oscstbounded}
Oscillation stable metric spaces are bounded.
\end{lem}
\begin{proof}
Let $\mathrm{M}=(M;\de)$ be an oscillation stable metric space and assume for a contradiction that $\mathrm{M}$ is unbounded. Fix a  sequence $(p_i; i\in \omega)$ of points in $M$ with:
\[
\forall (i\in \omega)\,  \big(\de(p_0,p_i)<\de(p_0,p_{i+1})\big) \text{  and} \lim_{i\to \infty}|\de(p_0,p_i)-\de(p_0,p_{i+1})|=\infty.
\]
Let $l(i)=\de(p_0,p_{i+1}-\de(p_0,p_i)$ and  $f: M\to \Re$ be  the function so that if   $\de(p_0,x)\in [\de(p_0,p_i),\de(p_0,p_{i+1})]$, then:
\begin{align*}
&f(x)= \frac{1}{l(i)}(\de(p_0,x)-\de(p_0,p_i)  \text{ if $i$ is even},\\
&f(x)=\frac{1}{l(i)}(\de(p_0,p_{i+1})-\de(p_0,p_x) \text{ if $i$ is odd}.
\end{align*}

The function $f$ is the composition of a uniformly continuous piecewise linear map from $\Re\to \Re$ with the uniformly continuous distance function $\de(p_0,x)$ and hence is uniformly continuous. The range of $f$ is a subset of  the interval $[0,1]$ with $f(p_i)=0$ if $i$ is even and $f(p_i)=1$ if $i$ is odd. 

Note: If $0\leq r\in \Re$ and $\epsilon>0$ then there exists an $i\in \omega$, even, so that $0\leq f(x)<\epsilon$ for all $x$ with $|\de(p_0,x)-\de(p_0,p_i)|\leq r$. If $0\leq r\in \Re$ and $\epsilon>0$ then there exists an $i\in \omega$, odd, so that $0\leq 1-f(x)<\epsilon$ for all $x$ with $|\de(p_0,x)-\de(p_0,p_i)|\leq r$.

Given $\epsilon>0$. Because $\mathrm{M}$ is oscillation stable there exists an isometric embedding $\alpha$ of $\mathrm{M}$ to a  copy $\mathrm{M}^\ast=(M^\ast;\de)$ of\/ $\mathrm{M}$ in $\mathrm{M}$ so that $|f(x)-f(y)|<\frac{1}{3}\epsilon$ for all $x,y\in M^\ast$.  Let $r=\de(p_0,\alpha(p_0))$ and $i$ even  be such that $0\leq f(x)<\frac{1}{3}\epsilon$ for all $x$ with  $|\de(p_0,x)-\de(p_0,p_i)|\leq r$. Then $0\leq f(\alpha(p_i))<\frac{1}{3}\epsilon$ because $|\de(p_0,\alpha(p_i))-\de(p_0, p_i)|=|\de(p_0,\alpha(p_i))-\de(\alpha(p_0), \alpha(p_i))|\leq r$. Let  $j$ be odd   such that $0\leq 1-f(x)<\frac{1}{3}\epsilon$ for all $x$ with  $|\de(p_0,x)-\de(p_0,p_j)|\leq r$. Then $0\leq 1-f(\alpha(p_j))<\frac{1}{3}\epsilon$ because $|\de(p_0,\alpha(p_j))-\de(p_0, p_j)|=|\de(p_0,\alpha(p_j))-\de(\alpha(p_0), \alpha(p_j))|\leq r$. 

But then we arrived at the contradiction $|f(\alpha(p_i))-f(\alpha(p_j))|>\frac{1}{3}\epsilon$.

\end{proof}

\begin{lem}
If for every partition of $M$ into two parts $(X,Y)$ and for every $\epsilon>0$ there exists a  copy $\mathrm{M}^\ast=(M^\ast,\de)$ of \/ $\mathrm{M}$ in $\mathrm{M}$ with $M^\ast\subseteq \big(X\big)_\epsilon$ or $M^\ast\subseteq \big(Y\big)_\epsilon$ then $\mathrm{M}$ is indivisible. 
\end{lem}
\begin{proof}
By induction on the number of parts in the partition. Let $(B_0,B_1,B_2,\dots, B_{n-1},B_n)$ be a partition of $M$ and $\epsilon>0$ be given. Let $X=\bigcup_{i\in n}B_i$ and $Y=B_n$.  If there exists a copy $\mathrm{M}^\ast=(M^\ast,\de)$ of $\mathrm{M}$ with $M\subseteq \big(B_n)_{\frac{\epsilon}{2}}$ we are done. Otherwise there exists a copy $\mathrm{M}^\ast=(M^\ast; \de)$ of $\mathrm{M}$ with $M^\ast\subseteq \big(X\big)_{\frac{\epsilon}{2}}$. 
\end{proof}

\begin{lem}\label{claim:unif1}
Let $\mathrm{M}=(M; \de)$ be a metric space and $(X,Y)$ a partition of $M$ and $f: M\to \Re$ with \[
f(x)=\inf\{\de(x,y)\mid y\in Y\}
\]
and for all $y\in Y$
\[
f(y)=\inf\{\de(x,y)\mid x\in X\}.
\]
then $f$ is uniformly continuous.
\end{lem}
\begin{proof}
Let $\epsilon>0$ and $\delta=\frac{1}{2}\epsilon$ and let $p,q\in M$ with $\de(p,q)\leq \frac{1}{2}\epsilon$. If $p\in X$ and $q\in Y$ then $f(p)\leq \frac{1}{2}\epsilon$ and $f(q)\leq \frac{1}{2}\epsilon$ and hence $|f(p)-f(q)|\leq \epsilon$. If $p,q\in X$ there exists a point $u\in Y$ with $\de(p,u)<f(p)+\frac{1}{2}\epsilon$ and hence $f(q)\leq \de(q,u)<f(p)+\epsilon$ implying $|f(q)-f(p)|<\epsilon$. 

\end{proof}

\begin{lem}\label{claim:unif2}
Let $\mathrm{M}=(M; \de)$ be a metric space and $k>0$ and $(X,Y)$ a partition of $M$ so that $\de(x,y)\geq \delta$ for all $x\in X$ and $y\in Y$. Then the function $f: M\to \Re$  with $f(x)=0$ for all $x\in X$ and $f(y)=1$ for all $y\in Y$ is uniformly continuous. 
\end{lem}
\begin{proof}
$|f(x)-f(y)|=0$ for all $x,y\in M$ with $\de(x,y)<\delta$. 

\end{proof}

\begin{thm}\label{thm:appequosc}
A metric space  $\mathrm{M}=(M; \de)$ is  oscillation stable if and only if it is   approximately indivisible.
\end{thm}
\begin{proof}
On account of Lemma \ref{claim:1} it remains to prove that if $\mathrm{M}$ is oscillation stable then it is approximately indivisible.

Let $\epsilon>0$ be given and $(X,Y)$ a partition of $M$ into two parts.

Let $f: M\to \Re$ be the function so that for all $x\in X$
\[
f(x)=\inf\{\de(x,y)\mid y\in Y\}
\]
and for all $y\in Y$
\[
f(y)=\inf\{\de(x,y)\mid x\in X\}.
\]
The function $f$ is uniformly continuous because of Lemma \ref{claim:unif1}.

Hence there exists a copy $\mathrm{M}^\ast=(M^\ast; \de)$ of $\mathrm{M}$ in $\mathrm{M}$ so that $|f(p)-f(q)|<\frac{\epsilon}{2}$ for all $p,q\in M^\ast$. If $f(p)<\epsilon$  for all $p\in X\cap M^\ast$ then $M^\ast\subseteq \big(Y\big)_\epsilon$. 

Let $p\in X\cap M^\ast $ with $f(p)\geq \epsilon$. If there is no $q\in Y\cap M^\ast$ then $M^\ast\subseteq X\subseteq \big(X\big)_\epsilon$. Let  $M^\ast\cap Y\not=\emptyset$. Let $x\in X\cap M^\ast$ and $y\in Y\cap M^\ast$ with $\de(x,y)<\frac{\epsilon}{2}$. Then $f(y)<\frac{\epsilon}{2}$ and hence $|f(p)-f(y)|>\frac{\epsilon}{2}$ a contradiction. It follows that $\de(x,y)\geq\frac{\epsilon}{2}$ for all $x\in X\cap M^\ast$ and all $y\in Y\cap M^\ast$. 

The function $g: M^\ast\to \Re$ with $g(x)=0$ for all $x\in X\cap M^\ast$ and $g(y)=1$ for all $y\in Y\cap M^\ast$ is uniformly continuous according to Lemma~\ref{claim:unif2}. It follows that there exists a copy $\mathrm{M}^{\ast\ast}=(\mathrm{M}^{\ast\ast}; \de)$ of $\mathrm{M}^\ast$ in $\mathrm{M}^\ast$ so that $|g(a)-g(b)|<\frac{1}{2}$ for all $a,b\in M^{\ast\ast}$. This in turn implies that $g(a)=g(b)$ for all $a,b\in M^{\ast\ast}$ and hence that $M^{\ast\ast}\subseteq X$ or $M^{\ast\ast}\subseteq Y$. 

\end{proof}

\section{Cantor sets and the 4-values condition}\label{sect:Cantor-4val}

Let $\mathcal{R}\subseteq \Re_{\geq 0}$. Theorem \ref{lthm:associative} states that if $\mathcal{R}$ is closed, then it satisfies the 4-values condition if and only if the operation $\oplus$ is associative on $\mathcal{R}$. If $\mathcal{R}$ is not a closed subset of $\Re_{\geq 0}$ then $\mathcal{R}$ can satisfy the 4-values condition but might not be closed under the operation $\oplus$. Let for example $\mathcal{R}=[0,1)\cup [2,3]$ then $\mathcal{R}$ satisfies the 4-values condition. We aim to show that there are quite interesting closed subsets $\mathcal{R}$ of the reals which satisfy the 4-values condition and hence quite intricate homogeneous metric spaces which are oscillation stable. Lemma \ref{Fact:1} together with Lemma \ref{lem:multiplying} are the basic tools to generate subsets of the reals which satisfy the 4-values condition. The sets $\mathcal{R}\subseteq \Re_{\geq 0}$ constructed below have 0 as a limit are closed and satisfy the 4-values condition. Hence, according to Theorem \ref{thm:characterization}, for each such set $\mathcal{R}$ there exists a unique Urysohn metric space $\UR$. 

The ordinary Cantor set obtained by successively removing the middle open third of the closed intervals does not satisfy the 4-values condition because  $(\frac{1}{3}\oplus_\mathcal{R}\frac{2}{9})\oplus_\mathcal{R}\frac{1}{9}=\frac{1}{3}$ but $\frac{1}{3}\oplus_\mathcal{R}(\frac{2}{9}\oplus_\mathcal{R}\frac{1}{9})=\frac{2}{3}$.  But if $(l_1,l_2,l_3,l_4,\dots)$ is a finite or infinite sequence of numbers with $\frac{1}{3}<l_i<1$ for all $i$ then the set obtained by successively removing the middle open intervals of length $l_i$ times the length of the interval, satisfies the 4-values condition. See below for a more precise account. 

\begin{lem}\label{Fact:1} Let $\mathcal{R}$ be closed and satisfy the 4-values condition and let $c> 0$, then:  The set $c\mathcal{R}=\{cr\mid r\in \mathcal{R}\}$ satisfies the 4-values condition.  The set $\{x\in \mathcal{R}\mid x\leq c\}$ satisfies the 4-values condition.  The interval $[0,1]$ satisfies the 4-values condition. 
\end{lem}
\begin{proof}
Let $c>0$ and $r,s\in \mathcal{R}$. Then $cr\oplus_{c\mathcal{R}}cs=c(r\oplus_\mathcal{R}s)$, implying the first assertion. Let $r,s,t\in \{x\in \mathcal{R}\mid x\leq c\}:=\mathcal{S}$ and let $m=\max(\mathcal{S})$. If $r+s\geq m$, then  $(r\oplus_\mathcal{S} s)\oplus_\mathcal{S}  t=m\oplus_\mathcal{S}  t=m$ and $s\oplus_\mathcal{S}  t\geq s$ and hence   $r\oplus_\mathcal{S} (s\oplus_\mathcal{S}  t)=m$. Similarly we obtain $(r\oplus_\mathcal{S} s)\oplus_\mathcal{S}  t=r\oplus_\mathcal{S} (s\oplus_\mathcal{S}  t)$ if $t+s\geq m$ or $r+t\geq m$. 

Let $r+s\leq m$ and $r+t\leq m$ and $s+t\leq m$. Then $r\oplus_\mathcal{S}s=r\oplus_\mathcal{R}s$ and $r\oplus_\mathcal{S}t=r\oplus_\mathcal{R}t$ and $t\oplus_\mathcal{S}s=t\oplus_\mathcal{R}s$.   Hence if  $(r\oplus_\mathcal{S}s)+t\geq m$ and $r+(s\oplus_\mathcal{S}t)<m$ then $(r\oplus_\mathcal{R}s)\oplus_\mathcal{R}t\geq m$ but $r\oplus_\mathcal{R}(s\oplus_\mathcal{R}t)<m$. Hence if $(r\oplus_\mathcal{S}s)+t\geq m$ then $r+(s\oplus_\mathcal{S}t)\geq m$ implying $(r\oplus_\mathcal{S}s)\oplus_\mathcal{S}t= m=r\oplus_\mathcal{S}(s\oplus_\mathcal{S}t)\geq m$. Similarly $(r\oplus_\mathcal{S}s)\oplus_\mathcal{S}t= m=r\oplus_\mathcal{S}(s\oplus_\mathcal{S}t)\geq m$ if $r+(s\oplus_\mathcal{S}t)\geq m$. If $(r\oplus_\mathcal{S}s)+t< m$ and $r+(s\oplus_\mathcal{S}t)<m$ then $(r\oplus_\mathcal{S}s)\oplus_\mathcal{S}t=(r\oplus_\mathcal{R}s)\oplus_\mathcal{R}t=r\oplus_\mathcal{R}(s\oplus_\mathcal{R}t)=r\oplus_\mathcal{S}(s\oplus_\mathcal{S}t)$.

The set $\Re_{\geq 0}$ satisfies the 4-values condition, because $r\oplus_{\Re_{\geq 0}} s=r+s$, implying that the set $[0,1]$ satisfies the 4-values condition.

\end{proof}

\begin{lem}\label{lem:multiplying}
If $\mathcal{R}\subseteq \Re_{\geq 0}$ is closed and satisfies the 4-values condition and $l>2\cdot\max(\mathcal{R})$ then the set $\{r+nl\mid n\in \omega \text{ and $r\in \mathcal{R}$}\}$ is closed and satisfies the 4-values condition. 
\end{lem}
\begin{proof}
Let $n,m\in \omega$ and $r,s\in \mathcal{R}$. Then $(nl+r)+(ml+s)=(n+m)l+r+s<(n+m)l+l$. Hence $(nl+r)\oplus(ml+s)=(n+m)l+(r\oplus s)$ and we have, for $n,m,k\in \omega$ and $r,s,t\in \mathcal{R}$:
\begin{align*}
&\big((nl+r)\oplus(ml+s)\big)\oplus(kl+t)=\big((n+m)l+ (r\oplus s)\big)\oplus (kl+t)=\\
&(n+m+k)l+\big((r\oplus s)\oplus t\big)=(n+m+k)l+\big(r\oplus (s\oplus t)\big)=\\
&(nl+r)\oplus\big((m+k)l+(s\oplus t)\big)=(nl+r)\oplus\big((ml+s)\oplus (kl+t)\big).
\end{align*}

\end{proof}

Given an interval $[a,b]$ of the reals and $0\leq w\leq 1$  let $[a\big( w \big) b]$ denote the subset of $[a,b]$ obtained by removing the middle open interval of length $ w (b-a)$ from $[a,b]$. Note that $[a\big(0\big)b]=[a,b]$ and  $[a\big(1\big)b]=\emptyset$ and that:
\begin{align}
&[a\big(w)b]=[a,\gamma(a,w,b)]\cup [\delta(a,w,b),b], \text{  for:}\\
&\gamma(a,w,b)=\frac{1}{2}\big((1+ w )a+(1- w )b\big),\\
&\delta(a,w,b)=\frac{1}{2}\big((1- w )a+(1+ w )b\big).
\end{align}
For finite sequences $\vec{w}=(w_i; i\in n\in \omega)$ we  define recursively the set of disjoint intervals $[a\big(\vec{w}\, \big)b]$.  If $\vec{w }$ is the empty sequence then $[a\big(\vec{ w }\, \big)b]=[a,b]$. If $\vec{w}=(w)$, the sequence consisting of a single entry, then $[a\big(\vec{w}\, \big)b]=[a\big(w\, \big)b]$. In general, for $\vec{w}_\ast=(w_i; 1\leq i\in n)$:
\[
[a\big(\vec{w}\, \big)b]=[a\big(\vec{w}_\ast\, \big)\gamma(a,l_0,b\, \big)] \cup[\delta(a,l_0,b),b]\big(\vec{w}_\ast\, \big)b].
\]
Note that the set $[0\big(\vec{w}\, \big)b]$ is a scaled version of the set $[0\big(\vec{w}\, \big)c]$, that is $[0\big(\vec{w}\, \big)b]=\frac{b}{c}\cdot[0\big(\vec{w}\, \big)c]$. Hence:

\begin{lem}\label{lem:fact2}
Let $\vec{w}=(w_i; i\in n\in \omega)$ be a sequence of numbers and $b,c$ be two positive numbers.  Then $[0\big(\vec{w}\, \big)b]$ satisfies the 4-values condition if and only if $[0\big(\vec{w}\, \big)c]$ satisfies the 4-values condition. 
\end{lem}
Note that the set $[a\big(\vec{w}\, \big)b+a]$ is a translation of the set $[0\big(\vec{w}\, \big)b]$, that is $[a\big(\vec{w}\, \big)b+a]=\{x+a\mid x\in [0\big(\vec{w}\, \big)b]\}$. Hence we obtain by induction on $n$:

\begin{lem}\label{lem:Cantor4val}
For all $n\in \omega$ and sequences  $\vec{w}=(w_i; i\in n\in \omega)$  with $\frac{1}{3}<w_i<1$ for all $i\in n$, the set $[0\big(\vec{w}\, \big)1]$ satisfies the 4-values condition. 
\end{lem}
\begin{proof}
The interval $[0,1]$ satisfies the 4-values condition. Let $\frac{1}{3}<w<1$ and $\vec{w}=(w_i; i\in n\in \omega)$.  If $[0\big(\vec{w}\, \big)1]$ satisfies the 4-values condition, the set $[0\big(\vec{w}\, \big)\frac{1}{2}(1-w)]$ satisfies the 4-values condition  and $[\frac{1}{2}(1+w)\big(\vec{w}\, \big)1]$ is an $\frac{1}{2}(1+w)$ translation of $[0\big(\vec{w}\, \big)\frac{1}{2}(1-w)]$. Hence $[0\big(\vec{w}\, \big)\frac{1}{2}(1-w)]\cup [\frac{1}{2}(1+w)\big(\vec{w}\, \big)1]$ satisfies the 4-values condition,    because $\frac{1}{2}(1+w)>2\cdot\frac{1}{2}(1-w)$ and therefore Lemma \ref{lem:multiplying} applies. 

Let $\vec{v}$ be the sequence $(w, w_0, w_1,w_2,\dots, w_{n-1})$. Then $[0,\big(\vec{v}\big)1]=  [0\big(\vec{w}\, \big)\frac{1}{2}(1-w)]\cup [\frac{1}{2}(1+w)\big(\vec{w}\, \big)1]$ and hence $[0\big(\vec{v}\big)1]$ satisfies the 4-values condition. 
\end{proof}

\begin{defin}\label{defin:Cantset}
Let $\vec{w}=(w_i;i\in \omega)$ be a sequence with $0<w_i<1$ for all $i\in \omega$. Then 
\[
[0\big(\vec{w}\big)1]=\bigcap_{n\in \omega}\vec{W}_n  \text{\  \  \    for $\vec{W}_n=(w_i\ i\in n)$}.
\]
\end{defin}

\begin{thm}\label{thm:Cantset4val}
Let $\vec{w}$ be a finite or infinite sequence with $\frac{1}{3}<w_i<1$ for  all indices $i$. Then the Cantor type set $[0\big(\vec{w}\big)1]$ satisfies the 4-values condition. 
\end{thm}
\begin{proof}
If $\vec{w}$ is finite, the Theorem follows from Lemma \ref{lem:Cantor4val}. Let $\vec{w}=(w_i;i\in \omega)$ and $\vec{W}_n=(w_i\ i\in n)$. Let $\mathcal{R}_n=[0\big(\vec{W}_n\big)1]$ and $\oplus_{\mathcal{R}_n}:=\oplus_n$ and let $\mathcal{R}_\infty=[0\big(\vec{w}\big)1]$ and $\oplus=\oplus_{\mathcal{R}_\infty}$. 

Note that if $x$ is a boundary point of $\mathcal{R}_n$ for some $n$,  then $x\in \mathcal{R}_m$  and it is a boundary point of $\mathcal{R}_m$ for all $m>n$ and hence an element of $\mathcal{R}_\infty$ and it is a boundary point of $\mathcal{R}_\infty$. This implies that for all $a,b\in \mathcal{R}_\infty$ exists an index $n$ so that $a\oplus_n b=a\oplus_m b$ for all $m\geq n$ and hence $a\oplus b=a\oplus_n b$. Therefore, for all triples $\{a,b,c\}$ of numbers in $\mathcal{R}_\infty$, there exists an index $n$ so that $\oplus_n$ agrees with $\oplus $ on the set $\{a,b,c\}$. Because $\oplus_n$ is associative:
\[
(a\oplus b)\oplus c=(a\oplus_n b)\oplus_n  c=a\oplus_n(b\oplus_n c)=a\oplus(b\oplus c). 
\]
\end{proof}
We did not use the full strength of Lemma \ref{lem:multiplying} when splitting closed intervals into two parts separated by an open interval.  We could have instead split the intervals into finitely many parts separated by open intervals of the same length. Hence for infinite sequences of the type $\vec{w}=\big((l_i,m_i);i\in \omega\big)$,  with $m_i$ giving the number of parts, yield Cantor type sets of the form $[0\big(\vec{w}\big)1]$.  Of course the construction can then also be extended to any countable sequence of this type, yielding quite intricate Urysohn type spaces.

\end{document}